
\magnification1200
\input amssym.def 
\input amssym.tex 
\def\SetAuthorHead#1{}
\def\SetTitleHead#1{}
\def\NoindentAfter{\everypar={\setbox0=\lastbox\everypar={}}}
\def\H#1\par#2\par{{\baselineskip=15pt\parindent=0pt\parskip=0pt
 \leftskip= 0pt plus.2\hsize\rightskip=0pt plus.2\hsize
 \bf#1\unskip\break\vskip 4pt\rm#2\unskip\break\hrule
 \vskip40pt plus4pt minus4pt}\NoindentAfter}
\def\HH#1\par{{\bigbreak\noindent\bf#1\medbreak}\NoindentAfter}
\def\HHH#1\par{{\bigbreak\noindent\bf#1\unskip.\kern.4em}}
\def\th#1\par{\medbreak\noindent{\bf#1\unskip.\kern.4em}\it}
\def\endth{\medbreak\rm}
\def\pf#1\par{\medbreak\noindent{\it#1\unskip.\kern.4em}}
\def\df#1\par{\medbreak\noindent{\it#1\unskip.\kern.4em}}
\def\enddf{\medbreak}
\let\rk\df\let\endrk\enddf
\let\Roster\bgroup\let\endRoster\egroup
\def\\{}\def\text#1{\hbox{\rm #1}}
\def\mop#1{\mathop{\rm\vphantom{x}#1}\nolimits}
\def\sam{}\def\prose{}\def\MaxReferenceTag#1{}
\def\qedbox{\vrule width2mm height2mm\hglue1mm\relax}
\def\qed{\ifmmode\qedbox\else\hglue5mm\unskip\hfill\qedbox\medbreak\fi\rm}
\let\SOverline\overline\let\Overline\overline
\def\Smallfonts{}

\let\Item\item
\def\cite#1{{\bf[#1]}}
\def\Em#1{{\it #1\/}}\let\em\Em
\def\Bib#1\par{\bigbreak\bgroup\centerline{#1}\medbreak\parindent30pt
 \parskip2pt\frenchspacing\par}
\def\endBib{\par\egroup}
\newdimen\Overhang
\def\rf#1{\par\noindent\hangafter1\hangindent=\parindent
     \setbox0=\hbox{[#1]}\Overhang\wd0\advance\Overhang.4em\relax
     \ifdim\Overhang>\hangindent\else\Overhang\hangindent\fi
     \hbox to \Overhang{\box0\hss}\ignorespaces}
\def\bbH{{\Bbb H}}
\def\bbN{{\Bbb N}}\def\bbQ{{\Bbb Q}}
\def\bbR{{\Bbb R}}\def\bbZ{{\Bbb Z}}
\def\Coordinates{\bigbreak\bgroup\parindent=0pt\obeylines}
\def\endCoordinates{\egroup}

\newcount\NRcount\NRcount0
\def\NR{\advance\NRcount1\Nr{\the\NRcount.}}
\def\IT{\advance\NRcount1\Item{\the\NRcount.}}

\def\SetRev{\def\track{fellow travel}
\def\tracks{fellow travels}
\def\tracked{fellow travelled}
\def\tracking{fellow traveller}}

\SetRev            
\def\implies{$\Rightarrow$}
\def\bull{\item{\raise1pt\hbox{\sam$\bullet$}}}
\def\Hyp{\bbH}
\def\cusps{\cal C}
\def\Hypbar{\Overline\bbH}

\def\CH{\mop{CH}}

\def\int{\mop{int}}
\def\Isom{\mop{Isom}}
\def\rank{\mop{rank}}

\def\CAT{\mop{CAT}}
\def\Autstruct{\mop{S\frak A}}
\def\AAutstruct{\mop{\frak A}} 
\def\BAutstruct{\mop{BS\frak A}}
\def\QAutstruct{\mop{Q\frak A}}
 
\def\len{\mop{len}}

\def\Par{{\cal P}}
\let\eval\Overline
\def\Abstract#1\par{\bgroup\Smallfonts\narrower\HHH #1\par}
\def\endAbstract{\par\egroup}
\def\TitleH{Geometrically Finite Hyperbolic Groups}
\def\Title{Automatic Structures, Rational Growth,\\ and \TitleH}
\def\\{}
\def\Author{Walter D. Neumann and Michael Shapiro}
\SetTitleHead{\TitleH}
\SetAuthorHead{\Author}

\H \Title
 
\Author\footnote{*}{Both authors
acknowledge support from the NSF for this research.}
 
\Abstract Abstract
 
We show that the set $\Autstruct(G)$ of equivalence classes of
synchronously automatic structures on a geometrically finite
hyperbolic group $G$ is dense in the product of the sets
$\Autstruct(P)$ over all maximal parabolic subgroups $P$. The set
$\BAutstruct(G)$ of equivalence classes of biautomatic structures on
$G$ is isomorphic to the product of the sets $\BAutstruct(P)$ over the
cusps (conjugacy classes of maximal parabolic subgroups) of $G$. Each
maximal parabolic $P$ is a virtually abelian group, so $\Autstruct(P)$
and $\BAutstruct(P)$ were computed in \cite{NS1}.

We show that any geometrically finite hyperbolic group has a
generating set for which the full language of geodesics for $G$ is
regular.  Moreover, the growth function of $G$ with respect to this
generating set is rational.  We also determine which automatic
structures on such a group are equivalent to geodesic ones. Not all
are, though all biautomatic structures are.
\endAbstract

\HH 1.~~Introduction

In \cite{NS1} the concept of equivalence of synchronous or
asynchronous automatic structures on a group (for definitions see
below) was introduced, and, among other things, the sets
$\Autstruct(G)$ and $\BAutstruct(G)$ of equivalence classes of
automatic or biautomatic structures on a group G were computed in
various situations. In this paper we describe the situation for
geometrically finite hyperbolic groups.  We also discuss the existence
of regular geodesic languages on such groups.

There are several definitions of what it means for a subgroup $G$ of
the isometry group $\Isom(\bbH^n)$ of hyperbolic $n$-space to be a
geometrically finite group. The traditional one is that there exist a
finite sided convex polyhedral fundamental domain for the action of
$G$ on $\bbH^n$.  Ratcliffe \cite{R} relaxes this to the requirement
that the convex polyhedron be locally finite-sided and each point $x$
in its closure in $\SOverline \bbH^n$ have a neighborhood which meets
only those faces $P$ incident to $x$.  It is unknown if these
definitions are equivalent; they are equivalent in dimensions $1$,
$2$, and $3$.  Epstein et al.\ in \cite{ECHLPT} take as their
definition that there exist an equivariant system of disjoint open
horoballs at the parabolic fixed points of $G$ and, moreover, if $X$
is the convex hull of the limit set with these horoballs removed then
$X/G$ is
compact. Ratcliffe's and Epstein's definitions are equivalent, as
follows from \cite{R, Theorem 12.6.2}.  Ratcliffe cites Bowditch
\cite{Bo} as an antecedent for this theorem.  We use the
Ratcliffe/Epstein concept of geometric finiteness.

Let $G$ be a geometrically finite hyperbolic group. Let $\Par$ be the
set of maximal parabolic subgroups $P\subset G$.  Each $P$ is a
finitely generated virtually free abelian group. Hence, as described
in \cite{NS1}, the set $\Autstruct(P)$ of equivalence classes of
automatic structures on $P$ is naturally equivalent to the set of
ordered rational linear triangulations of the sphere $S^{\rank(P)-1}$.
Given an element $[L]\in\Autstruct(G)$, we show there is an induced
element $[L_P]\in\Autstruct(P)$.  This induces a mapping
$$\Autstruct(G)\to\prod_{P\in\Par}\Autstruct(P).$$ We show this
mapping is injective with dense image. We also compute the set of
synchronous biautomatic structures on $G$: If $\cusps$ is a set of
representatives for the conjugacy classes of maximal parabolic
subgroups of $G$ then the mapping $$\BAutstruct(G) \to
\prod_{P\in\cusps} \BAutstruct(P)$$ is a bijection.

In \cite{NS2} we show an analogous result for $G$ equal to the
fundamental group of a graph of groups with finite edge groups.  In
that case the conjugates of the vertex groups play the role that the
maximal parabolics do in the geometrically finite hyperbolic case.
That result also holds for asynchronous automatic structures.

Similar results cannot hold for asynchronous automatic structures on
geometrically finite hyperbolic groups.  Indeed, in \cite{NS1, Sect.
4} we show that the set of asynchronous automatic structures on a
cocompact 3-dimensional hyperbolic group can be enormous, despite the
absence of cusps.  On the other hand, our results apply without change
to the set $\QAutstruct(G)$ of \em{quasigeodesic} asynchronous
automatic structures on $G$.  We suspect, but cannot prove, that this
set equals $\Autstruct(G)$ for a geometrically finite hyperbolic group
$G$.

In the early 80's Jim Cannon found examples of non-cocompact
geometrically finite groups with generating sets with rational growth
functions and asked if all such groups have such generating sets.
These examples, and his proof of rationality of the growth function
for cocompact hyperbolic groups (with respect to any generating set;
this holds for any word-hyperbolic group), were a major motivation in
the development of the theory of automatic groups. In \cite{ECHLPT} a
biautomatic structure is constructed for any geometrically finite
hyperbolic group.  If this structure consisted of word-geodesics it
would follow that the group has rational growth function for the given
generating set. The structure of \cite{ECHLPT} arises from a
``weighted geodesic structure'' on a groupoid containing the group but
this seems insufficient to guarantee rational growth function.

In Section 4 we show that any geometrically finite hyperbolic group
$G$ has a generating set $A$ so that the geodesics in $A^*$ form a
regular language and the growth function is rational (Theorem 4.3).
This is done by using a criterion which essentially goes back to
\cite{C}, namely, that any word in $A^*$ which is not geodesic has a
close neighbor which is shorter.  We call this criterion
``falsification by a fellow traveller.''

In Section 5, we describe which automatic structures on $G$ are
equivalent to automatic structures consisting of geodesics. This
depends on understanding the following problem.  Given a virtually
abelian group $P$, a generating set $A$, and an equivalence class of
automatic structures on $P$, when is there a geodesic structure
$L\subset A^*$ in this class?  The answer is encoded in a euclidean
polyhedron determined by translation lengths with respect to $A$.  It
turns out (Theorem 5.6) that not all automatic structures on
geometrically finite hyperbolic groups have geodesic representatives.
However, all the biautomatic ones do.

Our approach is to study $G$ by means of its action on the space $X$
described above.  Our viewpoint is in part inspired by the suggestive
but incomplete arguments of \cite{T}. Since that paper is unpublished,
we have collected our interpretation of those arguments in an appendix.

Several proofs in this paper could be greatly simplified if one had
affirmative answers to the following:

\th Questions

Suppose $G$ is a geometrically finite hyperbolic group.

{1.} Does there exist a finite index subgroup $H\subset G$, all of
whose parabolic subgroups are abelian?

{2.} If $P$ and $P'$ are distinct maximal parabolic subgroups of
$G$, does there exist a finite index subgroup $H\subset G$ in which
$P\cap H$ and $P'\cap H$ are non-conjugate?
\endth

This is essentially asking about the so-called ``LERF'' property for
maximal abelian parabolic subgroups respectively maximal parabolic
subgroups of $G$.

\HH 2.~~Background and definitions

Let $G$ be a finitely generated group and $A$ a finite set and
$a\mapsto\eval a$ a map of $A$ to a monoid generating set $\eval
A\subset G$.  As is usual, $A^*$ denotes the free monoid on $A$ and
the natural projection $A^*\to G$ is denoted $w\mapsto\eval w$.  Any
subset $L$ of $A^*$ which surjects onto $G$ is called a
\Em{normal form} for $G$. A \Em{rational structure} is a normal form
that is a regular language (i.e., the set of accepted words for some
finite state automaton, see below).

The \Em{Cayley graph} $\Gamma_A(G)$ is the directed graph with vertex
set $G$ and a directed edge from $g$ to $g\eval a$ for each $g\in G$
and $a\in A$; we give this edge a label $a$.  We do not require that
$\eval A=\eval A^{-1}$. It follows that the ``distance function''
$d(g,h)$ defined as the length of a shortest directed path from $g$ to
$h$ in $\Gamma_A(G)$ is not necessarily a symmetric distance function,
although it is positive definite and obeys the triangle inequality.
It exceeds the undirected distance from $g$ to $h$ in $\Gamma_A(G)$ by
at most a constant multiple; this constant is the size of the largest
$A$-word needed to express an element of $\eval A^{-1}$.  Thus
bounding directed distance is equivalent to bounding undirected
distance.  We denote $\ell(g)=d(1,g)$.

Each word $w \in A^*$ defines a path $[0,\infty)\to\Gamma$ in the
Cayley graph $\Gamma=\Gamma_A(G)$ as follows (we denote this path also
by $w$): $w(t)$ is the value of the $t$-th initial segment of $w$ for
$t=0,\ldots,\len(w)$, is on the edge from $w(s)$ to $w(s+1)$ for
$s<t<s+1\le\len(w)$ and equals $\eval w$ for $t\ge\len(w)$. We refer
to the translate by $g\in G$ of a path $w$ by $gw$.

A normal form $L$ for $G$ has the \Em{synchronous \tracking{}
property} if there exists a $\delta_L\in\Bbb N$ such that, given any
normal form words $v,w\in L$ with $\eval w = \eval{va}$ for some $a\in
A\cup A^{-1}\cup\{1\}$, the distance $d(w(t),v(t))$ never exceeds
$\delta_L$.  The normal form $L$ has the \Em{asynchronous \tracking{}
property} if $\delta_L\in\Bbb N$ exists such that for any $v$ and $w$ as
above there exists a non-decreasing proper function $t\mapsto t'\colon
[0,\infty)\to[0,\infty)$ such that $d({v(t)},{w(t')})\le\delta_L$ for
all $t$.  In these situations we say that $v$ and $w$
\Em{synchronously} (respectively \Em{asynchronously})
\Em{$\delta_L$-fellow travel}.

$L$ is a \Em{synchronous (asynchronous) automatic structure} for $G$
if it is a rational structure that has the synchronous (asynchronous)
\tracking{} property.  As discussed in \cite{NS1} following \cite{S1},
in the asynchronous case this definition is not quite equivalent to
the definition in \cite{ECHLPT} using automata, but it is equivalent
for finite-to-one languages, and every asynchronous automatic
structure $L$ (in either definition) contains a sublanguage which is a
one-to-one asynchronous automatic structure. 

``Automatic structure'' will mean ``synchronous automatic structure.''
An automatic structure $L$ is \Em{biautomatic} if its fellow traveller
constant $\delta_L$ can be chosen such that if $v,w\in L$ satisfy
$\eval w=\eval{av}$ with $a\in A$ then $av$ and $w$ synchronously
$\delta_L$-fellow travel.

Two asynchronous or synchronous automatic structures $L_1\subset
A_1^*$ and $L_2\subset A_2^*$ on $G$ are \Em{equivalent} if $L_1\cup
L_2 \subset(A_1\cup A_2)^*$ is an asynchronous automatic structure
(see \cite{NS1}).  Equivalently, there exists a $\delta$ such that
elements of $L_1$ and $L_2$ with the same value asynchronously
$\delta$-fellow travel.  We then write $L_1\sim L_2$.

Recall that a \Em{finite state automaton} $\cal A$ with alphabet $A$
is a finite directed graph on a vertex set $S$ (called the set of
\Em{states}) with each edge labelled by an element of $A$ and such
that different edges leaving a vertex always have different labels.
Moreover, a \Em{start state} $s_0\in S$ and a subset of \Em{accepted
states} $T \subset S$ are given.  A word $w\in A^*$ is in the
\Em{language $L$ accepted by $\cal A$} if and only if it defines a
path starting from $s_0$ and ending in an accept state in this graph.
We may assume there is no ``dead state'' in $S$ (a state not
accessible from $s_0$ or from which no accepted state is accessible).
Eliminating such states does not change the language $L$ accepted by
$\cal A$.
 
A language is \Em{regular} if it is accepted by some finite state
automaton.

Given a synchronous or asynchronous automatic structure $L$ on a group
$G$, a subset $N\subset G$ is \Em{$L$-rational} if $\{w\in L : \eval
w\in N\}$ is a regular sublanguage of $L$.  Proposition 2.7 of
\cite{NS1} says that if $L\sim L'$ then $N$ is $L$-rational if and
only if it is $L'$-rational. If $N$ is an $L$-rational subgroup of $G$
then there is a unique synchronous or asynchronous automatic structure
$L_N$ up to equivalence on $N$ such that $L$-words with value in $N$
fellow travel their corresponding $L_N$-words (cf.\ \cite{GS} or
\cite{NS1, Lemma 2.6}).

Given $\lambda\ge 1$ and $\epsilon\ge 0$, a map $f\colon X\to Y$ of
metric spaces is a \Em{$(\lambda,\epsilon)$-quasi-isometric map} if
$${1\over\lambda}d_Y(x,y)-\epsilon\le d_X(f(x),f(y))\le \lambda
d_Y(x,y)+\epsilon $$
for all $x,y\in X$.  If $X$ is an interval, we speak of a
\Em{quasigeodesic path} in $Y$.  Two metric spaces $X$ and $Y$ are
\Em{quasi-isometric} if there exists a quasi-isometric map $f\colon
X\to Y$ such that $Y$ is a bounded neighborhood of the image of $f$.
Then $f$ is called a \Em{quasi-isometry}.

\HH 3.~~Automatic structures on geometrically finite
hyperbolic groups

Let $G$ be a geometrically finite subgroup of $\Isom(\bbH^n)$.
We let $\Par$ be the
set of maximal parabolic subgroups $P< G$.  Each $P$ is a
finitely generated virtually abelian group.

Let $L$ be an automatic structure on $G$.  It is shown in \cite{S2}
that $P$ is $L$-rational. Thus, there is an induced element
$[L_P]\in\Autstruct(P)$ which only depends on $P$ and the class
$[L]\in\Autstruct(G)$.  We thus get a mapping
$$\Phi\colon\Autstruct(G)\to\prod_{P\in\Par}\Autstruct(P).$$ 

\th Theorem 3.1

$\Phi$ is injective with dense image.\endth

Before we proceed with the proof we must recall the geometry of the
situation, as described in \cite{ECHLPT}, following \cite{Bo}.  For a
detailed account, see \cite{R}.

Denote the standard compactification of $\Hyp^n$ by $\Hypbar^n$.
There is a smallest non-empty convex subset of $\Hyp^n$
on which $G$ acts. It is called the \Em{convex hull} for $G$ and
denoted $\CH(G)$. (It can be constructed as
$\CH(G)=C\Lambda(G)\cap\Hyp^n$, where $C\Lambda(G)$ is the convex hull
in $\Hypbar^n$ of the limit set $\Lambda(G)\subset\partial\Hypbar^n$.)

Each maximal parabolic subgroup $P$ fixes a unique point of
$\partial\Hypbar^n$. It preserves any horoball centered at this point.
We can choose a $G$-equivariant disjoint system of such horoballs, one
for each maximal parabolic subgroup.  Denote the horoball
corresponding to $P$ by $B_P$.  We will study $G$ via its action on $$
X = \CH(G) - \bigcup_{P\in\Par} \int(B_P) \prose.  $$ The quotient $M$
of $X$ by $G$ has the natural structure of a compact topological
orbifold with boundary.  In fact $\pi_1^{\mop{orb}}(M)=G$.  We call
the boundary piece $S_P=X\cap\partial B_P$ that resulted from removing
$\int(B_P)$ a \Em{horosphere of} $X$. In particular, there is a
one-one correspondence between maximal parabolic subgroups and
horospheres.  Two horospheres $S_P$ and $S_{P'}$ have the same image
in $M$ if and only if $P$ and $P'$ are conjugate in $G$.  It is known
that $G$ has finitely many conjugacy classes of maximal parabolic
subgroups $P$, so $M$ has finitely many horospherical boundary pieces.

As metric on $X$ we take the path metric, that is the metric given by
lengths of paths with length of each path computed using the standard
hyperbolic riemannian metric.  In particular, this makes $X$ a
geodesic metric space on which $G$ acts by isometries, cocompactly,
with finite stabilizers.

Let $A=A^{-1}$ be a finite generating set for
$G=\pi_1^{\mop{orb}}(M)$.  Choose a basepoint $\tilde p \in X$ with
trivial stabilizer.  Let $p\in M$ be its image.  Choose disjoint
closed paths at $p$ representing the elements of the generating set
$A$.  The inverse image in $X$ of the union of these loops is an
embedded copy $\Gamma$ of the Cayley graph $\Gamma_A(G)$.  It is a
standard result (cf.\ Milnor \cite{M}) that this embedding is a
quasi-isometry. In fact, the paths can be chosen so the embedding is a
$(\lambda,0)$-quasi-isometry for some $\lambda$. We identify the
Cayley graph with its image in $X$.

We now give a brief description of the geometry of the geodesics of
$X$ following Tatsuoka \cite{T}.  For more details see the appendix to
this paper.  The geodesics of $X$ consist of what Tatsuoka terms
``glancing geodesics''.  These are $C^1$-smooth paths made up of
hyperbolic geodesics off the horospheres of $X$ alternating with
euclidean geodesics lying on these horospheres.  Given such an
$X$-geodesic, $\gamma$, Tatsuoka studies the retraction $r_\gamma$ of
$X$ onto $\gamma$ mapping each point $x$ to its hyperbolically closest
point on $\gamma$.  In the appendix we take a somewhat different
viewpoint, retracting $X$ onto a neighborhood of $\gamma$ union those
horospheres that $\gamma$ meets.  We show that this is ``locally
strongly distance decreasing'' away from $\gamma$ and its horospheres
and deduce the following lemma (cf.\ Lemmas A5 and A6 of the
appendix).

\th Lemma 3.2

For any quasigeodesity constants $(\lambda,\epsilon)$ there exists a
constant $l$ such that if $w$ is a $(\lambda,\epsilon)$-quasigeodesic
in $X$ from $x$ to $y$ and $\gamma$ is the $X$-geodesic from $x$ to
$y$ then $w$ asynchronously $l$-fellow travels a path obtained from
$\gamma$ by possibly modifying $\gamma$ on those horospheres $S_P$
that $\gamma$ visits.  \qed
\endth

\pf Proof of injectivity in Theorem 3.1

Let $L\subset A^*$ and $L'\subset A^*$ be inequivalent finite-to-one
automatic structures on $G$.  We need to show that the induced
structures $L_P$ and $L'_P$ are inequivalent for some maximal
parabolic $P$.  Since the Cayley graph metric and $X$-metric on
$\Gamma$ are quasi-isometric, we may use the $X$-metric when
discussing the fellow traveller property.

There exist quasigeodesity constants $(\lambda,\epsilon)$ such that
the edge path in $\Gamma$ determined by any $w\in L\cup L'$ is a
$(\lambda,\epsilon)$-quasigeodesic in $X$.  Let $l$ be chosen by the
above Lemma.  We want to apply this Lemma for paths given by words
$w$.  We therefore increase this $l$ by the maximum $X$-length of a
generator so that, when we subdivide the path $w$ according to the
parts of $\gamma$ that it fellow travels, we can do so at letter
boundaries.  

Since $L$ and $L'$ are inequivalent, for each $k>2l$ we can find words
$x_k\in L$ and $x'_k\in L'$ so that $\eval{x_k} = \eval{x_k'}$ but
$x_k$ and $x_k'$ fail to asynchronously $k$-fellow travel with respect
to the $X$-metric. Let $\gamma_k$ be the $X$-geodesic from $1$ to
$\eval{x_k}$.  By Lemma 3.2, the failure of $x_k$ and $x'_k$
to fellow travel occurs in an $l$-neighborhood of some horosphere
$S_{P_k}$.  We can decompose $x_k$ and $x_k'$ as products of words
$$\eqalign{ x_k &= u_kv_kw_k\cr
x'_k &= u'_kv'_kw'_k\cr} $$
such that the portion $\eval{u_k}v_k$ of the path $x_k$ (that is, the
portion of $x_k$ labelled by $v_k$) begins and ends within distance
$l$ of the beginning and end of the portion of $\gamma_k$ lying on the
horosphere $S_{P_k}$, and similarly for the decomposition of $x'_k$.
This decomposition has the following properties:
\Roster

X1.~~the portion $\eval{u_k}v_k$ of the path $x_k$ and the portion 
$\eval{u'_k}v'_k$ of 
the path $x'_k$ both run in the $l$-neighborhood of $S_{P_k}$;

X2.~~these portions start and end at most $X$-distance $2l$ apart;

X3.~~these portions fail to asynchronously $k$-fellow travel.
\endRoster

For each $k$ we record the data
$\{s_k,s'_k,\eval{u_k}^{-1}\eval{u'_k}\}$, where $s_k$ is the state of
the machine for $L$ reached by $u_k$ and $s'_k$ is the state of the
machine for $L'$ reached by $u'_k$.  By property X2 above,
$\ell(\eval{u_k}^{-1}\eval{u'_k})\le\lambda 2l+\epsilon$. There are
therefore only finitely many possibilities for the data
$\{s_k,s'_k,\eval{u_k}^{-1}\eval{u'_k}\}$, so we can find infinitely
many values $k$ with the same data. By taking a subsequence and
renumbering, we may assume this data is the same for each $k$. In
particular, the equation $\eval{u_1}^{-1}\eval{u'_1}=
\eval{u_k}^{-1}\eval{u'_k}$ implies
$$
\eval{u_1}\eval{u_k}^{-1}=\eval{u'_1}\eval{u'_k}^{-1} \eqno{(*)}
$$
for all $k$.

We will now consider the paths
$$
y_k:=u_1v_k\quad\text{and}\quad y'_k:=u'_1v'_k\prose.
$$ 
Since $s_1=s_k$ we can choose a short word $t_k$ (of length at most
$|\cal A|$, where $\cal A$ is an automaton for $L$) such that
$y_kt_k\in L$. Similarly, we can find a short $t'_k$ such that
$y'_kt'_k\in L'$.  We claim that for each $k$ there exists a maximal
parabolic $Q_k$ such that:
 \Roster

Y1.~~the portion $\eval{u_1}v_k$ of the path $y_k$ and the portion 
$\eval{u'_1}v'_k$ of 
the path $y'_k$ both run in an $l$-neighborhood of $S_{Q_k}$;

Y2.~~these portions start and end at most $X$-distance $2l$ apart;

Y3.~~these portions fail to asynchronously $k$-fellow travel.
 \endRoster
Indeed, for Y1 note that the path $\eval{u_1}v_k$ runs in a
$l$-neighborhood of the horosphere $\eval{u_1}\eval{u_k}^{-1}S_{P_k}$,
which is the horosphere $S_{Q_k}$ with
$Q_k=\eval{u_1}\eval{u_k}^{-1}P_k\eval{u_k}\eval{u_1}^{-1}$.  Using
equation $(*)$, the same argument shows that the path
$\eval{u'_1}v'_k$ runs in a $l$-neighborhood of $S_{Q_k}$.  A similar
computation using $(*)$ deduces Y2 and Y3 from X2 and X3.

Now note that the horosphere $S_{Q_k}$ is at most $X$-distance
$\lambda(\ell(u_1)+\epsilon)+l$ from our basepoint in $X$, so there
are only finitely many possibilities for $Q_k$.  Thus, by taking a
subsequence once again, we may assume that $Q_k$ equals a fixed
maximal parabolic $P$ for all $k$. We claim that there is a uniform
bound, independent of $k$, on the distance of $\eval{y_k}$ and
$\eval{y'_k}$ from $P$.  This will complete the proof of injectivity,
since, under the assumption that the languages $L_P$ and $L'_P$ are
equivalent, $y_k$ and $y'_k$ would have to asynchronously
$\alpha$-fellow travel each other for some constant $\alpha$ depending
only on the above uniform bound and the fellow traveller constant
between $L_P$ and $L'_P$.

To see that $\eval{y_k}$ is a bounded distance from $P$, first note
that if $\Gamma_P$ is the graph consisting of edges of $\Gamma$ that
lie in the $l$-neighborhood of $S_P$, then the quotient
$P\backslash\Gamma_P$ is finite.  The points $\eval{u_1v_k}$ and
$\eval{u_1}$ of $\Gamma_P$ are connected in $\Gamma_P$ by the path
labelled by $v_k^{-1}$, so the image points in $P\backslash\Gamma_P$
can be connected in $P\backslash\Gamma_P$ by some path of length at
most $\mop{diam}(P\backslash\Gamma_P)$.  Let this path be labelled
$z_k$.  Then $\eval{u_1v_k}\eval{z_k}\in P\eval{u_1}$, so
$\eval{u_1v_k}\eval{z_k}\eval{u_1}^{-1}\in P$.  Thus
$\eval{y_k}=\eval{u_1v_k}$ is within $\Gamma$-distance
$\mop{diam}(P\backslash\Gamma_P)+\len(u_1)$ of $P$.  The same argument
applies to $\eval{y'_k}$.  \qed

We postpone the proof of dense image.  While this can be proven with
the tools now in hand, the technical details become much simpler once
we have developed some information about geodesic automatic
structures.

Now let $\BAutstruct(G)$ be the set of equivalence classes of
biautomatic structures on our geometrically finite hyperbolic group
$G$.  Given $[L]\in \BAutstruct(G)$, the induced structure $L_P$ on a
maximal parabolic subgroup $P$ is biautomatic.  Moreover, since
biautomaticity implies invariance under conjugation, $[L_P]$
determines $[L_Q]$ for each conjugate $Q$ of $P$.  Theorem 3.1 thus
leads to an injective map
$$\BAutstruct(G) \to
\prod_{P\in\cusps} \BAutstruct(P),$$
where $\cusps$ is a set of representatives for the conjugacy classes
of maximal parabolic subgroups of $G$.  On the other hand,
\cite{ECHLPT} shows that a biautomatic structure can be constructed on
$G$ from any choice of one biautomatic structure at each cusp, in
other words, the above map is surjective (see also Remark 5.7).  Thus:

\th Theorem 3.6

The above map $\BAutstruct(G) \to\prod_{P\in\cusps} \BAutstruct(P)$
is a bijection.\qed\endth

The proof of Theorem 3.1 goes through with no essential change if
$\Autstruct(G)$ is replaced by the set $\QAutstruct(G)$ of equivalence
classes of quasigeodesic asynchronously automatic structures on
$G$. In fact, we suspect, but cannot prove, that these two sets are
equal.

\HH 4.~~Language of geodesics

Let $A$ be a monoid generating set for a group $G$.  A word $w\in A^*$
is \Em{geodesic} if it has shortest length among $A$-words
representing $\eval w$.  We shall use the following result to detect
when the language of geodesic words is regular.

\th Proposition 4.1

If $G$ is a group with finite monoid generating set $A$ and there
exists a $\delta$ such that any non-geodesic directed path in the
Cayley graph has a shorter directed path with the same value that
asynchronously $\delta$-fellow travels it, then the language of
geodesic words is regular.
\endth

\pf Proof

The proof recalls the standard proof that the language of geodesics in
a word hyperbolic group is regular (\cite{ECHLPT}) which has its
origins in \cite{C}.  To test if a path $u$ is geodesic, as we move
along $u$ we must keep track at each time $t$ of what points $x$ in a
$\delta$-neighborhood of $u(t)$ have been reached by paths that
asynchronously $\delta$-fellow travel $u$.  We must also record the
optimal time differential to reach the point $x$ by such a path.  The
point $x$ can clearly be reached in time at most $t+k\delta$, by
following $u$ to $u(t)$ followed by a geodesic path to $x$ (here $k$
is the constant relating $d(g,h)$ to $d(h,g)$).  If $x$ is reached in
time less than $t-\delta$ then $u$ is clearly not geodesic.  Thus, the
relevant time differential lies in the interval $\{-\delta,
-\delta+1,\ldots,k\delta-1,k\delta\}$ and the information that must be
kept track of is the function $\phi\colon B(\delta)\to\{-\delta,
-\delta+1,\ldots,k\delta-1,k\delta\}$, where $B(\delta)$ is the ball
of radius $\delta$ in the Cayley graph.  We can build a finite state
automaton with the set of such maps as states plus one ``fail state''.
The initial state is the map $\ell$ (recall $\ell(x)=d(1,x)$).  The
$a$-transition from a state $\phi$ leads to the following state $\psi$
if this $\psi$ satisfies $\psi(1)=0$ and to the fail state otherwise.
$$
\psi(x)=\Bigl\{
 \eqalign{&\phi(ax)-1 \cr &\min\{\phi(y)| y\in B(\delta), d(y,ax)=1\}\cr}
 \eqalign{&\text{~~if}\quad ax\in B(\delta)\prose,\cr
 &\text{~~if}\quad ax\notin B(\delta)\prose,\cr}\Bigr.
$$
for $x\in B(\delta)$.  \qed

\df Definition

We say that a monoid generating set $A$ for a group $G$ that satisfies
the premise of the above proposition has the \Em{falsification by
fellow traveller property}.
\enddf

\rk Question

Can one find a monoid generating set $A$ of a group $G$ so the
language of geodesic words is regular but $A$ does not have the
falsification by fellow traveller property?
\endrk

\th Proposition 4.2

If $A$ has the falsification by fellow traveller property then the
growth function of $G$ with respect to $A$ is rational.\endth

\pf Proof

Recall that the growth function in question is the power series
$$
f(t)=\sum_{g\in G}t^{\ell(g)}\prose.
$$
It is \Em{rational} if it is the power series expansion of a rational
function of $t$.

It is a standard result that the growth of a regular language $L$ is
rational.  If $\cal A$ is a finite state automaton for $L$ one forms
the transition matrix $M$ for $\cal A$ whose rows and columns are
indexed by the states of $\cal A$ and whose entry $m_{ij}$ counts the
number of edges from state $i$ to state $j$. Then the number of words
of length $n$ in $L$ is $v_1M^nv_2$ where $v_1$ is the row vector with
a $1$ at the start state and $0$'s elsewhere and $v_2$ is the column
vector with $1$'s at accept states and $0$'s elsewhere.  The growth
function is then given by the rational function
$v_1\bigl(\sum_{i=0}^\infty (tM)^i\bigr)v_2=v_1(I-tM)^{-1}v_2$ (see,
e.g., \cite{C}).

If $L$ is a geodesic language for $G$ which does not biject to $G$,
the growth of $L$ clearly overcounts the growth of $G$. We can
compensate for this overcount in the following way.

Let $\cal A$ be the machine constructed in the previous proof and
$L$ the language of geodesics accepted by this machine.  We will call
$g'\in G$ a ``parent'' of $g\in G$ if there exists an outbound edge in
$\Gamma$ from $g'$ to $g$.  We claim that the number of parents of
$\eval w$, $w\in L$, is determined by the state of $\cal A$ reached by
$w$.  For if $\psi$ is the state, this number is the number of $h\in
B(\delta)$ such that $\psi(h)=-1$ and there is a directed edge from
$h$ to $1$ (we must assume here that $\delta$ has been chosen at least
as large as $k$).  We can correct the overcount in the previous
paragraph by replacing the matrix $M$ by $M'$ with entries
$m'_{ij}=m_{ij}/p_j$, where $p_j$ is the number of parents of state
$j$. We can assign $p_j$ arbitrarily at the start and fail states,
since this does not affect $v_1(M')^nv_2$.   The growth
function is thus the rational function $v_1(I-tM')^{-1}v_2$.
\qed

\th Theorem 4.3

If $G$ is a geometrically finite hyperbolic group then $G$ has a
generating set $A=A^{-1}$ with the falsification by fellow traveller
property.  In particular, the set of geodesic words forms a regular
language and the growth function is rational.
\endth

\rk Remark

The conclusion of this Theorem will not be true for every
generating set.  In fact, J.~Cannon has given the following example of
a generating set for a virtually abelian group such that the geodesic
language is not regular.  Consider the split extension $P$ of
$\bbZ^2$, generated by $\{a,b\}$, by $\bbZ/2$, generated by $t$, such
that $t$ conjugates $a$ to $b$ and $b$ to $a$.  As generators of $P$
we shall take $a,a^{-1},c,c^{-1},d,d^{-1},t,t^{-1}$ with $\eval
c=a^2,\eval d=ab$.  Then a word of the form $tc^ntc^m$ is geodesic so
long as $m<n$, but can be replaced by the shorter word $d^{2n}c^{m-n}$
if $m\ge n$.  But it is easy to see that a regular language $L$ that
contains $tc^ntc^{n-1}$ must also include words of the form $tc^ntc^m$
with $m>n$ if $n-1$ exceeds the number of states of a machine for $L$.

Since the maximal parabolics in a geometrically finite hyperbolic
group $G$ are virtually abelian, we must take care to avoid this sort
of behavior.\endrk

\th Proposition 4.4 

Any finite generating set for an abelian group has the falsification
by fellow traveller property. If $P$ is a virtually abelian group
then any finite generating set for $P$ is contained in a
generating set $A$ with the falsification by fellow traveller
property. $A$ may be chosen with $A=A^{-1}$.
\endth

\pf Proof

We first consider the case that $P$ is abelian. Let $A=\{a_1,\ldots,
a_m\}$ be a monoid generating set.  If $u\in A^*$ we will denote the
total exponent of $a_i$ in $u$ by $n_i(u)$ and ${\bf
n}(u):=(n_1(u),\ldots,n_m(u))$.  We write
$(n_1,\ldots,n_m)\le(n'_1,\ldots,n'_m)$ if $n_i\le n'_i$ for each
$i=1,\ldots,m$. We claim that there exists a bound $k$ such that if
$u$ is a non-geodesic word then one can find ${\bf
n}=(n_1,\ldots,n_m)\le {\bf n}(u)$ with $\sum_in_i\le k$ for which
$a^{\bf n}:=a_1^{n_1}\dots a_m^{n_m}$ is non-geodesic.  This does what
is required, for if $u_0$ is obtained from $u$ by deleting $n_i$
instances of the letter $a_i$ for each $i$ and $u_1$ is a geodesic
word with value $a^{\bf n}$ then $\eval u =\eval{u_0u_1}$ and $u_0u_1$
$2k$-fellow travels $u$.

To see the claim we first observe that $\bbN^m$ with this ordering has
the property that any subset has only finitely many minimal elements.
For if not, we would have an infinite sequence of pairwise
non-comparable elements in the lattice of $m$-tuples of natural
numbers.  But this cannot exist, since this lattice has the property
that any infinite sequence ${\bf n}_j$ in it has a subsequence $\{{\bf
n}_{j_l}\}$ with ${\bf n}_{j_{l}}\le{\bf n}_{j_{l+1}}$ for all $l$.
Indeed, it is clear that the lattice of natural numbers has this
property, and that a finite product of lattices with this property has
this property.  Now let $S$ be the set of $m$-tuples giving rise to
non-geodesic words.  Take $k$ to be the maximal coordinate sum of a
minimal element of $S$.  This proves the claim.

Now suppose $P$ is given by a short exact sequence 
$$
0\to N\to P\to F\to 1
$$ 
with $N$ abelian and $F$ finite. Let $B$ be a monoid generating set.
By enlarging $B$ if necessary we may assume that $\eval B$ surjects
onto $F$ under $P\to F$.  For any $w\in B^*$ we may then find $n_w\in
N$ such that $\eval w = n_w$ if $\eval w\in N$ and $\eval w\in
n_w\eval B$ otherwise.  Let $C$ be a generating set of $N$ which
includes $n_w$ for every $w$ of length at most 3, contains
$\eval{b{b'}^{-1}}$ for any $b,b'\in B$ with $\eval{b{b'}^{-1}}\in N$,
and is mapped into itself by all inner automorphisms of $P$.  We claim
that the generating set $A=B\cup C$ has the falsification by
fellow traveller property.  Indeed, if $u\in A^*$ is a
non-geodesic word which has no $B$-letters in it then we have already
shown that it is fellow travelled by a shorter word.  If it has at
least three $B$-letters in it, then we use the invariance of $C$ under
inner automorphisms to move the last three $B$-letters in $u$ to the
end of $u$. They then form a three-letter terminal segment $w$, which
we can replace by $n_wb$ with $b\in B$ to obtain a word $v$ with
$\len(v)=\len(u)-1$ and $\eval v=\eval u$.  This word $v$
fellow travels $u$.  By the same argument, if $u\in A^*$ has
one or two $B$-letters in it then $u$ is fellow travelled by a word
$vb\in C^*B$ of the same length and value.  If $\eval b\in N$ then we
can again apply the abelian case already proved.  If not, there is a
geodesic $v'b'\in C^*B$ from $1$ to $\eval u$.  Note that
$\len(v')<\len(v)$. Since $\eval{b{b'}^{-1}}\in N$ we can find $c\in
C$ such that $\eval{cb'}=\eval{b}$.  Then
$\eval{vcb'}=\eval{vb}=\eval{v'b'}$, so $\eval{vc}=\eval{v'}$. By
replacing $vc$ at most twice by a shorter $C$-word with the same value
that fellow travels it we replace $vcb'$ by a word which
fellow travels $u$, is shorter than $u$, and has value $\eval u$.\qed

We shall need some preparation for the proof of Theorem 4.3.

\th Lemma 4.5

Let $G$ be a group and $A$ a generating set with the falsification by
fellow traveller property. Then if $u$ is a geodesic and $g,h\in G$,
then there exists a geodesic $v$ with value $\eval v= g\eval u h$ such
that the paths $u$ and $g^{-1}v$ asynchronously fellow travel with
fellow traveller constant $2(\delta_A+1)(\ell(g)+\ell(h))$, where
$\delta_A$ is the constant guaranteeing the falsification by fellow
traveller property.  \endth

\pf Proof

It is not hard to give a synchronous version of this, but we will not
need this.

Suppose $g=1$ and $h=\eval a$ with $a\in A$.  Then a geodesic with
value $\eval{ua}$ has length at least $\len(ua)-2$, so $ua$ can be
turned into a geodesic $v$ by at most two repeats of replacing it by a
shorter path that $\delta$-fellow travels it.  Thus the Lemma is
proved in this case.  If $g=\eval a$ with $a\in A$ and $h=1$ then the
same argument applied to $au$ proves the Lemma.  The general case is
now an induction on $\ell(g)$ and $\ell(h)$.\qed

Now let $B$ be a generating set for $G$.  Let $P_1,\ldots, P_m$ be a
set of representatives for the conjugacy classes of maximal parabolic
subgroups of $G$.  Given any constant $K$, Proposition 4.4 implies
that we may for each $i$ choose a generating set $A_i$ for $P_i$ with
the falsification by fellow traveller property and containing any
element of $P_i$ which moves the basepoint of $X$ at most distance
$K$.  Let $A=B\cup\bigcup_{i=1}^m A_i$. We may include in $A_i$ any
elements of $B$ which happen to evaluate into $P_i$.  Since distinct
parabolic subgroups are disjoint, $A_i$ is then the set of $a\in A$
which evaluate into $P_i$.  We claim that, if $K$ is large enough, the
set $A$ is a generating set for $G$ with the falsification by fellow
traveller property.

Let $\Gamma$ be the Cayley graph for $G$ with respect to $A$, embedded
in $X$ as in the previous section.  

\th Lemma 4.6

Suppose $K$ is large enough and $A$ is as above. Then, for any $l$
there is a $K_l$ such that, if $S_P$ is any translate of the
horosphere $S_{P_i}$ and $w$ is a $\Gamma$-geodesic segment which
travels entirely in an $l$-neighborhood of $S_P$, then $w$ is labelled
by a word $u_0pu_1$ with $p\in (A_i)^*$ and $\len(u_0),\len(u_1)\le
K_l$.
\endth

\pf Proof 

By performing a translation we may assume $P=P_i$.

First note that there is a retraction $\rho\colon X\to S_{P}$ along
$\bbH^n$-geodesics perpendicular to $S_{P}$. This follows from the
fact that the convex hull $\CH(G)$ is convex; these geodesics are the
geodesics pointing to the parabolic fixed point of $P$.

Choose $k$ large enough that the image of any $B$-edge from the
basepoint of $X$ lies in the $k$-neighborhood $N_k(S_{P_j})$ of each
of the horospheres $S_{P_j}$, $j=1,\ldots,m$.  Consider the above
retraction $\rho\colon X\to S_{P}$. There is an overall bound $D$ on
the diameter of the $\rho$-images of the $N_k(S_Q)$ in $S_{P}$ as $Q$
runs through the maximal parabolic subgroups of $G$ other than $P$.
To see this, we use the upper half space model and put the fixed point
of $P$ at $\infty$, and make $S_{P}$ the horizontal hyperplane which
lies at height $1$ in the model.  For each $Q\ne P$, $S_Q$ lies below
this hyperplane, so $N_k(S_Q)$ lies below a hyperplane distance $k$
above $S_{P}$.  Since the projection $\rho$ is by vertical lines, the
claim is now immediate.

By our choice of $k$, if $e$ is an edge that does not connect two
points of $P$, it lies in $N_k(S_Q)$ for some $Q\ne P$, so its
$\rho$-image has length at most $D$.  

Let $d=\max_{x\in S_{P}}\min_{p\in P} d_X(x,p)$ (recall we are
identifying $\Gamma$, and hence $G$, with a subset of $X$).  We shall
take $K$ sufficiently large that any element of $P=P_i$ that moves the
basepoint at most $3D+2d$ is in $A_i$.  Let $(\lambda,\epsilon)$ be
the quasigeodesity constants relating Cayley graph distance and
$X$-distance.

Suppose $l$ is chosen.  We will show that
$K_l=3\lambda(l+d)+3\epsilon+2$ satisfies the lemma.

Suppose first that $w$ is a geodesic path in $N_l(S_{P})$ which does
not meet $P$.  We must show it has length less than $K_l$.  Consider
its projection $\rho w$ onto $S_{P}$, which has length at most
$D\len(w)$ since each letter of $w$ moves at most distance $D$ in the
projection.  Let $x_0,\ldots,x_m$ be points spaced at most $3D$ apart
along the path $\rho w$ from the beginning point $x_0$ to the end
point $x_m$.  We can take $m\le D\len(w)/(3D) + 1=\len(w)/3 + 1$.  For
each $x_t$ let $y_t$ be a point of $P$ within distance $d$ of it.
Then the successive $y_t$'s differ by elements of $A_i$.  Let $v$ be
the path so determined.  We can get from each endpoint of $w$ to the
corresponding endpoint of $v$ by a $\Gamma$-path of length at most
$\lambda(l+d)+\epsilon$.  Thus, we have constructed a path from the
beginning point of $w$ to its end point of length at most
$2\lambda(l+d)+2\epsilon+ \len(w)/3  +1$.  Since $w$ was geodesic,
$\len(w)\le 2\lambda(l+d)+2\epsilon + \len(w)/3+1$, which implies
$\len(w)\le3\lambda(l+d)+3\epsilon+3/2<K_l$.

Now suppose $w$ is a geodesic path in $N_l(S_{P})$ which starts and
ends in $P$.  We will show it never strays from $P$.  Suppose to the
contrary that $w$ is a shortest counterexample.  We apply the argument
of the preceding paragraph. The term $2\lambda(l+d)+2\epsilon$ now does
not appear, so we get the inequality $\len(w)\le \len(w)/3 +1$ giving
$\len(w)\le 3/2$.  This implies $\eval w\in A_i$, which is a
contradiction.  This completes the proof of the lemma.\qed

\th Lemma 4.7  

Choose the generating set $A$ for $G$ as in Lemma 4.6.  Then there
exist a $\delta$ and $l$ such that a path $u$ which has no shorter
path that $\delta$-fellow travels it satisfies the conclusion of Lemma
3.2. That is, if $\gamma$ is the $X$-geodesic from the initial
point of $u$ to the end-point of $u$, then $u$ asynchronously
$l$-fellow travels a path obtained from $\gamma$ by possibly modifying
$\gamma$ inside $l$-neighborhoods of those horospheres $S_P$ that
$\gamma$ visits.
\endth

We postpone the proof of this Lemma and first show how Theorem 4.3
follows from it.  

\pf Proof of Theorem 4.3 

We will show that $G$ has the falsification by fellow traveller
property.  Given  $l$ and $\delta$ as in Lemma 4.7, any larger $l$ and
$\delta$ also work.  We choose such an $l$ and $\delta$, but may
increase $\delta$ later. 

Lemma 4.6 now implies that our path $u$ can be written
in the form $u_0p_0u_1p_1\ldots$ such that each subword $p_j$ is a
(possibly empty) word in the elements of the parabolic generating set
$A_{i_j}$ corresponding to the $j$-th horosphere that the 
$X$-geodesic $\gamma$ visits. Moreover, the path $u$ fellow travels
$\gamma$, except possibly along these parabolic portions
$\eval{u_0\ldots u_j}p_j$ of the path $u$.  Since $u$ is not
$\delta$-fellow travelled by a shorter path, by assuming $\delta$ is
larger than the $\delta$'s for the $P_i$'s we ensure that the portions
$p_j$ of $u$ must be geodesic.

If $v$ is a geodesic path with value $\eval v=\eval u$ then $v$ has a
similar decomposition $v_0q_0v_1q_1\ldots$, and each portion
$\eval{v_0\ldots v_j}q_j$ of $v$ begins and ends a bounded distance
from the beginning and end of the portion $\eval{u_0\ldots u_j}p_j$ of
$u$.  Since $\eval{u_0\ldots u_j}p_j$ is geodesic, Lemma 4.5 lets us
replace this portion $\eval{v_0\ldots v_j}q_j$ of $v$ by a new
parabolic geodesic which fellow travels $\eval{u_0\ldots u_j}p_j$ at
some appropriate distance determined by the Lemma 4.5 and the bounds
that have occurred so far.  We may assume $\delta$ is larger than this
bound.  This replacement does not change the length of $v$.  Doing
this for each $j=0,1,\ldots$ replaces $v$ by a geodesic which
$\delta$-fellow travels $u$.  Since $u$ had no shorter fellow
traveller, it must itself be geodesic, and we are done.\qed

\pf Proof of Lemma 4.7 

Given a path of length at most $2\delta$, any geodesic with the same
endpoints lies in a $\delta$ neighborhood of it.  Thus the assumption
on $u$ implies that any $2\delta$-long subpath of $u$ is geodesic.  We
will modify $u$ by replacing maximal horospherical segments of $u$ by
their corresponding hyperbolic geodesics.  We will show that resulting
path $\hat u$ (which no longer lives in $X$) is a local quasigeodesic
in $\bbH^n$.  Appeal to Lemma 4.8 will show that this is a (global)
quasigeodesic.  It therefore fellow travels its $\bbH^n$-geodesic
$\gamma$.  Lemma 4.7 will then follow.  We now provide the details.

We take $(\lambda,0)$ to be the quasi-isometry constants
relating $\Gamma$ and $X$.  We suppose $\delta$ is larger than the
falsification by fellow traveller constants for the parabolic
subgroups $P_i$.  Then any subword of $u$ which lies in any $A_i^*$ is
geodesic.  For each $i$ we replace every maximal $A_i^*$ substring
with the corresponding $\bbH^n $-geodesic, and call the resulting path
in $\bbH^n$ $\hat u$.  We claim that there is $(M, q)$ so that every
$\log(2\delta/\lambda)$-long subpath of $\hat u$ is an $\bbH^n$
$(M,q)$-quasigeodesic.  Let $\mu$ be such a subsegment of $\hat u$.

First suppose that both ends of $\mu$ lie in $X$.  Then $\mu=\hat v$
for some subpath $v$ of $u$.  We have, say, $v=v_1a_1v_2 \ldots
a_{j-1}v_j$, where the $a_i$ are the horospherical segments of $v$.
Thus $\mu=\hat v = v_1\mu_1v_2 \ldots \mu_{j-1} v_j$, where the
$\mu_i$ are the corresponding hyperbolic geodesics.  Using Lemma A3
which compares hyperbolic and $X$-distance
we then have
$$\len(\mu) \ge 2\log(\len(v)/\lambda)>\log(\len(v)/\lambda).$$ 
This forces $\len(v)\le 2 \delta$ and thus $v$ is a Cayley graph
geodesic.  We suppose that $\nu$ is the $\Hyp^n$-geodesic for $v$.
Then by Lemmas 3.2 and A6, $\mu$ travels in a bounded neighborhood of
$\nu$.  We now check that there are global quasigeodesity constants
$(\lambda',\epsilon')$ so that each $v_i$ is an $\Hyp^n$
$(\lambda',\epsilon')$-quasigeodesic.  Since each of these is a Cayley
graph geodesic, it is an $X$ $(\lambda,0)$-quasigeodesic.  Let
$\nu_i$ be the $\Hyp^n$-geodesic for any subpath $v'_i$ of $v_i$.
Notice that $\nu_i$ cannot stray more than a bounded amount into any
horoball, for otherwise a long portion of $v'_i$ would lie on a
horosphere and thus have been replaced in $\mu=\hat v$.  Thus there is
a bound on the ratio between the $X$-distance between endpoints of
$v'_i$ and the $\Hyp^n$-distance between these endpoints.  This gives
us the desired constant $\lambda'$.

We now observe that $\mu$ consists of $\Hyp^n$-geodesics and
$\Hyp^n$-quasigeodesics all of which travel in a bounded corridor of
$\nu$.  By an argument similar to Cannon's ``progression in geodesic
corridors'' \cite{C}, this makes $\mu$ an $\Hyp^n$
$(M,q_1)$-quasigeodesic.  $M$ and $q_1$ depend only on the fact that
$\delta$ exceeds the falsification by fellow traveller constants of
the parabolics in our given generating set.  Accordingly, we can
increase $\delta$ without changing $M$ and $q_1$.

We now turn to the case where one or more ends of $\mu$ penetrate a
horoball, but only do so by a bounded amount.  That is, we have $\mu =
\mu_1 u_1 \ldots u_{j-1} \mu_j$, where each $u_i$ is a subword of $u$,
each $\mu_i$ is a hyperbolic geodesic inside a horoball, and $\mu_1$
and $\mu_j$ are of length at most $q'$.  (One of these may be empty.)
We choose $q'$ so that any hyperbolic segment of length greater than
$q'$ which contacts the horosphere and stays within the horoball
makes an angle of close to $\pi/2$ with the horosphere.  It is an easy
exercise to see that if one affixes a path of length at most $q'$ to
an $(M,q'')$-quasigeodesic, the resulting path is an
$(M,q''')$-quasigeodesic, where $q'''=q''+(M+1)q'$.  Since
$u_1\ldots u_{j-1}$ is an $(M,q_1)$-quasigeodesic, it follows that
$\mu$ is an $(M,q_2)$-quasigeodesic where $q_2=q_1+(2M+2)q'$.

Finally we must address the case in which one or both of the ends of
$\mu$ lies inside a horoball and is long, i.e., $\mu = \mu_1 u_1
\ldots u_{j-1} \mu_j$ and one or both of $\mu_1$ and $\mu_{j-1}$ has
length at least $q'$.  So suppose $\len(\mu_1) \ge q'$.  Then by
choice of $q'$, $\mu_1$ makes an angle of close to $\pi/2$ with its
horosphere.  Let $\sigma$ be the $\Hyp^n$-geodesic for $u_1 \dots
u_{j-1}$.  By our first case, $u_1 \dots u_{j-1}$ is an
$(M,q_1)$-quasigeodesic, and thus stays close to $\sigma$.  It follows
that $\sigma$ cannot stray far into the horoball that $\mu_1$
penetrates, for otherwise a long initial segment of $u_1$ would lie
close the horosphere.  By Lemma 4.6, $u_1$ would include parabolic
generators and thus would not appear in $\mu$.  Thus $\sigma$ either
lies outside the horoball of $\mu_1$, or enters it at an angle which
is bounded away from $\pi/2$.  Thus the angle between $\mu_1$ and
$\sigma$ is bounded away from $0$ by some positive constant $\alpha$.
Thus, there is a constant $q_\alpha$ so that $\mu_1\sigma$ a
$(1,q_\alpha)$-quasigeodesic.  Accordingly, $\mu_1 u_1 \dots u_{j-1}$
is an $(M,q_3)$-quasigeodesic where $q_3$ depends only on $M$, $q_2$,
and $q_\alpha$.

We now suppose that $\mu_j$ is also long.  If $u_1 \dots u_{j-1}$ is
sufficiently long, then standard results of hyperbolic geometry show
that $\mu_1^{-1}$ and $\mu_j$ emanating from the endpoints of $u_1
\dots u_{j-1}$ and must diverge from each other and $\mu$ is an $(M,
q_4)$-quasigeodesic, where $q_4$ depends only on the previous
constants and the length of time necessary for each of $\mu_1^{-1}$
and $\mu_j$ to diverge from $u_1 \dots u_{j-1}$.  On the other hand,
if $u_1 \dots u_{j-1}$ is short, then $\mu_1$ and $\mu_j$ enter nearby
horoballs from nearby points at angles near $\pi/2$ and thus also
diverge from each other.  In this case, $\mu$ is a
$(1,q_5)$-quasigeodesic, where $q_5$ depends only on the bound for the
length of $u_1 \dots u_{j-1}$ and the angle (necessarily close to
$\pi/2$) that $\mu_1$ and $\mu_j$ make with their respective
horospheres.

It now follows that any $\log(2\delta/\lambda)$ subpath of $\hat u$ is
an $(M,q)$-quasigeodesic in $\Hyp^n$, where $q=\max \{q_1, q_2, q_3,
q_4, q_5\}$.  As we remarked above, we are free to increase $\delta$
without changing $M$ and $q$.  By taking $\delta$ sufficiently large,
it follows by Lemma 4.8 below that $\hat u$ is an
$(M_0,q_0)$-quasigeodesic in $\Hyp^n$.  In particular, $\hat u$ lies
close to its hyperbolic geodesic.  Lemma 4.7 now follows.\qed

We say a path $\sigma$ is a $k$-local $(\lambda,
\epsilon)$-quasigeodesic if every subpath of $\sigma$ of length at
most $k$ is a $(\lambda,\epsilon)$-quasigeodesic.

\th Lemma 4.8

Suppose $Y$ is a $\delta$-hyperbolic space. Given quasigeodesity
constants $(\lambda,\epsilon)$, there are $k$ and quasigeodesity
constants $(\lambda',\epsilon')$ so that every $k$-local
$(\lambda,\epsilon)$-quasigeodesic in $Y$ is a
$(\lambda',\epsilon')$-quasigeodesic in $Y$.
\endth

\pf Proof

Here we will use a ``parametrized'' version of $\delta$-hyperbolic
metric spaces.  Thus, if $\alpha\beta\gamma$ is a geodesic triangle,
then $\alpha$, $\beta$ and $\gamma$ decompose as
$\alpha=\alpha_0\alpha_1^{-1}$,$\beta=\beta_0\beta_1^{-1}$,and
$\gamma=\gamma_0\gamma_1^{-1}$, so that
$\len(\gamma_1)=\len(\alpha_0)$, $\len(\alpha_1)=\len(\beta_0)$, and
$\len(\beta_1)=\len(\gamma_0)$, and each of the pairs $\gamma_1,
\alpha_0$, $\alpha_1,\beta_0$ and $\beta_1,\gamma_0$ synchronously
$\delta$-fellow travel.  (See, e.g., \cite{ABC+}.)  

Recall that there is $\epsilon''$ so that every
$(\lambda,\epsilon)$-quasigeodesic stays within $\epsilon''$ of its
geodesic.  Suppose now, that $v$ and $v'$ are
$(\lambda,\epsilon)$-quasigeodesics emanating from a common point $p$
and that $v^{-1}v'$ is also a $(\lambda,\epsilon)$-quasigeodesic.  Let
$\nu$ and $\nu'$ be geodesics for $v$ and $v'$.  It now follows that
$\nu$ and $\nu'$ can $2\delta$-fellow-travel for distance at most
$t=\lambda\delta+\lambda\epsilon''+\epsilon/2 +\epsilon''$.
Indeed, let $q$ and $q'$ be points distance $t_0$ from $p$ along $\nu$
and $\nu'$ and suppose they are within distance $2\delta$ of each
other.  Let $r$ and $r'$ be points of $v$ and $v'$ distance at most
$\epsilon''$ from $q$ and $q'$ respectively.  Consider the portion $B$
of $\nu^{-1}\nu'$ from $r$ to $r'$.  Its endpoints lie at most
$2\delta+2\epsilon''$ apart, so
$\len(B)\le\lambda(2\delta+2\epsilon'')+\epsilon$.  On the other hand,
$\len(B)\ge2(t_0-\epsilon'')$.  These inequalities imply $t_0\le
\lambda\delta+\lambda\epsilon''+\epsilon/2 +\epsilon''=t$, as claimed.

We let $k$ be an even integer with $k/2\ge\lambda(2t+1)+\epsilon$.

We suppose $u$ is a subpath of a $k$-local $(\lambda,
\epsilon)$-quasigeodesic.  We first suppose that $\len(u)$ is a
multiple of $k/2$, and write $u=u_1\dots u_m$ with $\len(u_i)=k/2$ for
$i=1,\ldots,m$.  For each $u_i$, let $\mu_i$ be the corresponding
geodesic, which therefore has length $\ge2t+1$.  Let $\alpha_i$ be the
geodesic from the beginning of $u_1$ to the endpoint of $u_i$.  We
will show inductively that $\len (\alpha_i) \ge i$, and that
$\alpha_i^{-1}$ and $\mu_i^{-1}$ $\delta$-fellow travel for distance
at least $\len (\mu_i)-t\ge t+1$.  It will then follow that the
endpoints of $u$ are separated by distance at least $m=2\len (u)/k$.
Thus this $u$ is a $(2/k,0)$-quasigeodesic.  It then follows that even
if $\len(u)$ is not a multiple of $k/2$, then $u$ is a $(\lambda',
\epsilon')$-quasigeodesic in $Y$, where $\lambda'=2/k$ and
$\epsilon'=1+2/k$.

Our inductive hypotheses hold for $i=1$, so we must prove the
inductive step. 

Consider the triangle $\alpha_i\mu_{i+1}\alpha_{i+1}^{-1}$.
Notice that $\mu_{i+1}$ and $\alpha_i^{-1}$ cannot $\delta$-fellow travel
for more than distance $t$, for otherwise $\mu_{i+1}$ and $\mu_i^{-1}$
$2\delta$-fellow travel for this distance, contradicting our
observation about $\nu$ and $\nu'$ above.  Consequently
$\mu_{i+1}^{-1}$ and $\alpha_{i+1}^{-1}$ $\delta$-fellow travel for at
least distance $\len (\mu_{i+1})-t\ge t+1$.  Likewise $\alpha_{i}$ and
$\alpha_{i+1}$ $\delta$-fellow travel for distance at least
$\len (\alpha_{i})-t$, and this forces $\len (\alpha_{i+1})\ge i+1$.
This completes the induction.\qed 

The interested reader can check that by increasing $k$, we can force
$\lambda'$ as close as we like to $\lambda$.

\HH 5.~~Geodesic automatic structures

Let $G$ be a geometrically finite hyperbolic group.  In this section
we shall investigate which classes in $\Autstruct(G)$ can be
represented by geodesic languages (if $A$ is a monoid generating set
then $L\subset A^*$ is \Em{geodesic} if it consists of geodesic
words).  The following lemma will let us look for ``almost geodesic''
languages instead. $L\subset A^*$ is an \Em{almost geodesic language}
if there is a bound $K$ such that each $u\in L$ is at most $K$ longer
than a geodesic word representing $\eval u$.

\th Lemma 5.1

For any group, if $A$ is a monoid generating set with the
falsification by fellow traveller property and $L\subset A^*$ is an
almost geodesic automatic structure on $G$ then there exists a
geodesic automatic structure $L'\subset A^*$ with $L'\sim L$.\endth

\pf Proof

By applying the falsification by fellow traveller property at most $K$
times to a word $u$ of $L$ we may replace it by a geodesic word $w$
with $\eval w=\eval u$ which $K\delta$-fellow travels $u$.  The
language
$$
\{(u,w) : u\in L, \text{$w$ geodesic, $\eval u=\eval w$, $u$ and $w$
$K\delta$-fellow travel}\}
$$
is clearly the language of an asynchronous two-tape automaton (the
argument here recalls the standard comparator automata of
\cite{ECHLPT}).  Projection on the second factor is thus a regular
language (cf.\ e.g., \cite{S1}) and is the language $L'$ we
desire. \qed

We will first need to discuss geodesic automatic structures for a
finitely generated virtually abelian group $P$.  Such a $P$ is given
by an exact sequence
$$
1\to \bbZ^m\to P\to F\to 1
$$
with $F$ finite.  We will need to think of $\bbZ^m$ as being a subset
of $\bbR^m$.  For this reason we will often write the group structure
in $\bbZ^m$ additively.

Recall from \cite{NS1} that an automatic structure
$[L]\in\Autstruct(P)$ determines (and is determined by) a rational
ordered triangulation $\cal T_{L}$ of the sphere $S^{m-1}$ of linear
rays from the origin in $\bbR^m$.  The set of vertices of this
triangulation is denoted $\partial L$ and consists of the rays in $\bbR^m$
that are fellow travelled by rays of $L$ (a \Em{ray} of $L$ is an
infinite word, all of whose initial segments are initial segments of
$L$-words).

\th Proposition 5.2

Suppose $[L] \in \Autstruct(P)$.  Then the following are equivalent.
 \Roster \Item{1.}$[L]$ has a geodesic
representative $L'\subset A^*$ for some generating set $A$ of $P$;

\Item{2.}$[L]$ has an almost geodesic representative
$L'\subset A^*$ for some generating set $A$ of $P$;

\Item{3.} There is an $F$-invariant subset $S\subset \partial L$ which
lies in no hemisphere of $S^{m-1}$.
\endRoster
 Moreover, with the extra restriction $A=A^{-1}$, {\rm1} and
{\rm2} are equivalent and equivalent to:
 \Roster
\Item{3\sam${}'$.}There is an $F$-invariant subset $S\subset \partial
L$ such that $S\cap-S$ lies in no hemisphere of $S^{m-1}$.
\endRoster
\endth 

To prove this proposition we will need some preparation.  Let $A$ be a
finite monoid generating set for $P$.  For $v\in P$ let $\ell(v)$ be
the shortest length of an $A$-word representing $v$.  We write the
group structure in $\bbZ^m$ additively.  For $v\in\bbZ^m$ let
$\tau(v)=\lim_{n\to\infty}(\ell(nv)/n)$.  This is the \Em{translation
length} of $v$ as defined in \cite{GS}.  Trivially,
$\tau(cv)=c\tau(v)$ for $c\in\bbZ^+$. It follows that $\tau$ extends
to $\bbQ^m$ by $\tau(cv)=c\tau(v)$ for rational $c\ge0$.  By
sub-additivity of translation length, $\tau$ is continuous on
$\bbQ^m$, so we may extend it by continuity to $\bbR^m$.  Let
$$C(A):=\{v\in\bbR^m:\tau(v)\le1\}.$$
The notation is chosen to suggest ``convex hull'', in view of the
following lemma.

\th Lemma 5.3

If $F=\{1\}$, so $P=\bbZ^m$, then $C(A)$ is the convex hull of $A$
considered as a subset of $\bbR^m$. In general, $C(A)\subset \bbR^m$
is a rational polyhedron (polyhedron with rational vertices) with $0$
in its interior. It is invariant under the action of $F$ on $\bbR^m$.
\endth

\pf Proof

We will actually prove a more general version of the lemma. Suppose
each element $a$ of our generating set $A$ is assigned a positive
integral weight $\len(a)$. The length $\len(w)$ of a word is then
defined as the sum of weights of letters of $w$ and $\ell(g)$ is then
again defined as the shortest length of a word representing $g$.  The
situation of Lemma 5.3 is that all elements of $A$ have weight 1.  We
will show that the lemma holds for any weights.  In particular, in the
free abelian case $F=1$ we will show that $C(A)$ is the convex hull of $$
V(A):=\{{\scriptstyle{1\over\len(a)}}\eval a : a\in A\}\subset
\bbR^m\prose.\eqno(*) $$

We first note that $C(A)$ is the closure of $C(A)\cap\bbQ^m$. Indeed,
near any point $x\in C(A)$ we can find a rational point $y$ with
$\tau(y)$ close to $\tau(x)$ and hence $\tau(y)\le1+\epsilon$ for some
small $\epsilon$.  By multiplying $y$ by a rational number just below
$1/(1+\epsilon)$, we replace it by a rational point that is still
close to $x$ and is in $C(A)$.

We start with the special case that $F=1$, so $P=\bbZ^m$.  By the
above comment, we need only verify that rational points of $C(A)$ and
the convex hull of $V(A)$ agree.  Certainly, $V(A)$ is in $C(A)$.
Since translation length $\tau$ is sub-additive, it follows that the
convex hull of $V(A)$ is in $C(A)$. Conversely, suppose $x\in
C(A)\cap\bbQ^m$. Then for any $\epsilon>0$ we can find an integer
$n>0$ such that $nx\in\bbZ^m$ and $nx=\sum\lambda_i\eval a_i$ with
$(1/n)\sum\lambda_i\len(a_i)\le 1+\epsilon$. In particular, passing to
the limit as $n\to\infty$ gives $x=\sum\mu_i\eval a_i$ with
$\sum\mu_i\len(a_i)\le 1$. So $x$ is in the convex hull as claimed,
proving the lemma in the case that $P$ is free abelian.

Now suppose $P$ is not free abelian.  We will introduce an expanded
weighted generating set with the same translation function $\tau$.
Let $f=|F|$. For each word $w\in A^*$ of unweighted length $\le f$
which evaluates into $\bbZ^m$ we add a new element $a_w$ to $A$ with
weight $\len(w)$ and value $\eval a_w=\eval w$.  Denote this new
generating set by $A'$. Replacing $A$ by $A'$ does not alter weighted
geodesic length, so it does not alter $\tau$.  We now expand our
generating set again by adding a generator $h(a)$ for each $a\in A'$
with $\eval a\in \bbZ^m$ and each $h\in F-\{1\}$; we put
$\eval{h(a)}=h(\eval a)$ and $\len(h(a))=\len(a)$.  Denote this new
generating set by $A''$.  These new generators may change weighted
geodesic length, but we claim they do so by a bounded amount, so
translation length remains unchanged.

Indeed, geodesic length is certainly not increased, so we must just
show it is at worst decreased by a bounded amount.  So suppose $g\in
P$ is expressed by an $A''$-geodesic $w$.  Write $A''=A_Z\cup A_F$,
where $A_Z$ consists of the letters that evaluate into $\bbZ^m$ and
$A_F$ consists of the remaining ones.  Note that $A_F\subset A$.
Using the $F$-invariance of $A_Z$, we may move all $A_F$-letters to
the end of $w$.  Then, if the terminal segment of $w$ consisting of
$A_F$-letters has length $\ge f$, it has a subsegment which evaluates
into $\bbZ^m$. We can replace this by a letter of $A_Z$ and move it to
the front of $w$.  Repeating eventually gives a word consisting of
$A_Z$ letters followed by at most $f$ $A_F$ letters.  Moreover, if
$g\in\bbZ^m$ then there are no $A_F$ letters.  Each $A_Z$-letter is of
the form $h(a)$ for some $h\in F$ ($h(a)$ will mean $a$ for $h=1$).
Since they commute, we may collect together all letters of the form
$h(a)$ for each given $h$.  That is, our $A''$-geodesic now has the
form $w=w_1\ldots w_f u$, where $w_i$ is a word in the letters
$h_i(a)$ and $u$ is a word of $A_F^*$.  We can rewrite this in terms
of $A'$ as $w'=u_1w'_1u_1'\ldots u_fw'_fu_f' u$, where $u_i$ and
$u_i'^{-1}$ evaluate into the coset $f_i\bbZ^m$ and $w'_i$ is the word
in $(A')^*$ from which $w_i$ was made by $f_i$.  Replacing $w$ by $w'$
increases length by at most a bounded amount, to whit
$\sum_i\len(u_i)+\len(u_i')$, as claimed.

For $g\in\bbZ^m$, we have seen that translation length $\tau(g)$ is
the same whether computed using $A$ or $A''$.  But if computed using
$A''$, we have seen that only elements of $A_Z$ are needed, so
$\tau(g)$ is translation length with respect to $A_Z$.  We have thus
reduced to the free abelian case where the result is already proved. \qed

\th Lemma 5.4

For any $F$-invariant rational polyhedron $Q\subset \bbR^m$ containing
$0$ in its interior we can find a generating set $A$ of $P$ with the
falsification by fellow traveller property such that $C(A)=N.Q$ for
some $N>0$. (Here $N.Q:=\{Nx:x\in Q\}$.)  If $Q=-Q$ we may choose this
$A$ with $A=A^{-1}$.
\endth

\pf Proof

Let $B$ be a generating set of $P$. Construct $A=B\cup C$ with the
falsification by fellow traveller property, as in the proof of Lemma
4.4, with $C\subset\bbZ^m$. Let $A_Z$ be the subset of $A$ evaluating
into $\bbZ^m$.  As described in the proof of 4.4, any geodesic word
representing an element of $\bbZ^m$ lies in $A_Z^*$.  It follows that
$C(A)=C(A_Z)$.  Thus $C(A)$ is just the convex hull of $A_Z\subset
\bbR^m$. In the proof of 4.4 the only requirement on $C$ was that it
be $F$-invariant and contain certain elements, but any larger
$F$-invariant subset $C\subset\bbZ^m$ also works. We can thus replace
$C$, and hence $A_Z$, by any larger $F$-invariant subset of $\bbZ^m$.

Choose a large integer $N$ so that the vertices of $N.Q$ have integer
coordinates and $C(A_Z)\subset N.Q$. Replacing $C$ by $N.Q$ does what
is required.  If $Q=-Q$ and we started with a set $B$ with $B=B^{-1}$,
then the resulting $A$ has $A=A^{-1}$.\qed

\df Definition

The generating set $A$ constructed in the above proof has the property
that if $A_Z$ is the subset of $A$ evaluating into $\bbZ^m$ then $A_Z$
is $F$-invariant and every element of $\bbZ^m$ has an $A$-geodesic
representative involving only $A_Z$-letters.  We shall call such a
generating set \Em{good}.\enddf

Let $\partial A$ denote the points of $S^{m-1}$ represented by the
rays in $\bbR^m$ through the vertices of $C(A)$.

\th Lemma 5.5

Let $[L]\in \AAutstruct(P)$ and let $A$ be a monoid generating set of
$P$.  If there is a representative $L\in[L]$ which is an
almost-geodesic sublanguage of $A^*$ then $\partial A \subset
\partial L$. The converse holds if $A$ is good. 
\endth

\pf Proof

Suppose $L$ is an almost geodesic automatic structure on $P$.  We may
assume that $L$ bijects to $G$ by extracting a sublanguage if
necessary. In \cite{NS1} it is shown that any bijective automatic
structure on a virtually abelian group is a finite union of languages
of the form $N=\{u_0\}\{w_1\}^*\{u_1\}\ldots \{w_k\}^*\{u_k\}$, with
$u_j,w_j\in A^*$ and such that the rays defined by $\eval
w_1,\ldots,\eval w_k$ define an ordered simplex of the triangulation
$\cal T_L$.  (Finite unions of languages of this type are often called
 ``simply starred''; they are precisely the regular languages of
polynomial growth.) Moreover, there exists such an $N$ for each
maximal simplex of $\cal T_L$. Since $N$ is an almost geodesic
language, the words $w_1,\ldots, w_k$ must be geodesic.  Consider an
element $g=\eval w_1^{n_1}\ldots \eval w_k^{n_k}$.  For any $n$ the
value of the word $v_n:=u_0w_1^{nn_1}u_1\ldots w_k^{nn_k}u_k\in N$ is
close to $g^n$, so it fellow travels any $L$-representative of $g^n$,
so $\ell(g^n)$ differs from $\len(v_n)$ by a bounded amount.  Hence
$\ell(g^n)$ differs from $\sum nn_i\len(w_i)$ by a bounded amount, so
$\tau(g)=\sum n_i\len(w_i)$. It follows that the simplex in $\bbR^n$
spanned by ${1\over \len(w_i)}\eval w_i$, $i=1,\ldots,k$ lies on the
boundary of $C(A)$.  We have thus shown that every maximal simplex of
$\cal T_L$ corresponds to a face or portion of a face of $\partial
C(A)$, so the vertices of $C(A)$ can only occur at vertices of $\cal
T_L$, as claimed.

Now suppose $A$ is a good generating set and suppose all vertices of
$C(A)$ lie along rays determined by points of $\partial L$. Since $A$
is good we know that $C(A)=C(A_Z)$, where $A_Z$ is the subset of $A$
evaluating into $\bbZ^m$. Thus, $C(A)$ is the convex hull of the set
$V(A_Z)$ defined in the proof of Lemma 5.3.

For any $x\in \bbZ^m$, the ray through $x$ passes through a face of
$C(A)$, and some multiple $nx$ is then an integral linear combination
of the elements $a\in A_Z$ that determine this face.  This gives a
geodesic representative for $nx$ of length $\tau(nx)$.  Thus, by
taking a positive integral multiple $N.C(A)$ of $C(A)\subset \bbR^m$,
we can ensure that each ray determined by a point of $\partial L$
intersects $\partial (N.C(A))$ in a point $v$ of $\bbZ^m$ which has a
geodesic representative $w_v\in A_Z^*$ of length $\tau(v)$.  If
$\sigma=\langle [v_1],\ldots,[v_m]\rangle$ is an ordered simplex of
the triangulation $\cal T_{[L]}$ we let
$$
L_\sigma=\{w_{v_1}^{n_1}\ldots w_{v_m}^{n_m} : n_i\ge 0\}.
$$
We take $L''=\bigcup_\sigma L_\sigma$. This $L''$ consists of
geodesics, and its image in $\bbZ^m$ contains a finite index subgroup $H$
(for each $\sigma$ take the subgroup generated by the $\eval w_{v_i}$
and then intersect these).  Let $X\subset A^*$ be a finite set such
that $\eval X$ is a set of cosets representatives for $H$ in $P$.
Then $L'=L''X$ is an automatic structure in $[L]$ by $\cite{NS1}$, and
it clearly consists of almost geodesics.\qed

\pf Proof of Proposition 5.2

We first show {1$\Leftrightarrow$2}. Trivially {1\implies2}.
Conversely, {2} implies {1} by Lemma 5.1 if we can assume our
generating set $A$ has the falsification by fellow traveller property.
But we can assume this: by Lemma 5.4 we can replace $A$ by a good
generating set with the falsification by fellow traveller property at
the expense of multiplying $C(A)$ by some integer, and Lemma 5.5
implies that the new $A$ will still satisfy {2}.

Note that the image of $L$ under conjugation by $f$ is $fLf^{-1}$
which is equivalent to $fL$.  Thus, $\partial(fL)$ is the image of
$\partial L$ under conjugation by $f$ and the set $S:=\bigcap_{f\in
F}\partial({fL})$ is the maximal $F$-invariant subset of $\partial L$.

Now if {2} holds then Lemma 5.5 plus the $F$-invariance of $\partial
A$ implies that $\partial A\subset S$.  Since $C(A)$ has $0$ in its
interior and it is the convex hull of its vertices, its vertices
cannot lie in a half-space of $\bbR^m$.  Thus $\partial A$ cannot lie
in a hemisphere of $S^{m-1}$, so the same holds for $S$.  Conversely,
if $S$ does not lie in a hemisphere and if we choose rational points
on the rays in $\bbR^m$ defined by the points of $S$, the convex hull
$Q$ of these points will be a polyhedron containing $0$ in its
interior.  We can do this $F$-equivariantly.  Lemma 5.4 then gives a
good generating set $A$ with $C(A)=N.Q$ for some $N$.  Since $\partial
A$ is given by the vertices of $Q$, we have $\partial A\subset
S\subset\partial L$, so {2} holds by Lemma 5.5.

The corresponding statements under the restriction $A=A^{-1}$ follow
easily.\qed

Now suppose $G$ is a geometrically finite hyperbolic group. Recall
that if $L$ is an automatic structure on $G$ then we have induced
structures up to equivalence on each maximal parabolic subgroup $P$.
We denoted these structures $L_P$. Given a conjugate $Q=gPg^{-1}$ of
$P$, the language $g^{-1}L_Qg$ is a language on $P$ which we will
denote $L_P^g$. It is not hard to see, using Section 2, that there are
just finitely many different languages $L_P^g$ up to equivalence
(for each $P$ they number at most the number of states in a machine
for $L$).  $L$
is equivalent to a biautomatic structure if and only if for each $P$
the languages $L_P^g$, $g\in G$, are all equivalent to each other
(cf.\ Theorem 3.6).

Recall that $\cusps$ denotes a set of conjugacy representatives for
the maximal parabolic subgroups. The following is the main theorem of
this section.

\th Theorem 5.6

Suppose $G$ is a geometrically finite hyperbolic group and $[L] \in
\Autstruct(G)$.  Then the following are equivalent: \Roster \Item{1.} $G$
has a geodesic automatic structure equivalent to $L$; \Item{2.} $G$
has an almost geodesic automatic structure equivalent to $L$;
\Item{3.} for each $P\in\cusps$ the set $\bigcap_{g\in G}\partial
L_P^g$ is not contained in a hemisphere.  \endRoster In particular,
these hold if $L$ is biautomatic.  \endth

\rk Remark

L. Reeves \cite{Re} has used the above result to show that a subgroup
of a geometrically finite hyperbolic group $G$ is again geometrically
finite if and only if it is rational for some biautomatic structure on
$G$.
\endrk

\pf Proof

Clearly 1 implies 2.  We now show that 2 implies 3.  By going to a
sublanguage if necessary we will assume $L$ bijects to $G$.

We start by observing that we have two notions of translation length.
There is the geodesic translation length $\tau=\tau_A$ used above:
$\tau_A(g)=\lim_{n\to\infty}\ell(g^n)/n$.  There is also the language
translation length $\tau_L(g)=\liminf_{n\to \infty}\len(w_n)/n$, where
$w_n\in L$ is chosen with $\eval{w_n}=g^n$ for each $n$.

Note that $\tau_A$ is a conjugation invariant function on $G$.  On the
other hand, if $L$ is almost geodesic, then $\tau_A=\tau_L$, so
$\tau_L$ is also conjugation invariant.  

Now fix a maximal parabolic subgroup $P\subset G$ and let $1\to\bbZ^m
\to P\to F\to 1$ give its structure.  The translation length $\tau_L$
induces a translation length on $\bbZ^m\subset P$.  We extend this to
a map $\tau:\bbR^m\to\bbR_+$ and define $C(P)=\{x:\tau(x)\le
1\}\subset \bbR^m$ as before.  Denote $L(P)=\{w\in L:\eval w\in P\}$.
Then, as in the proof of Lemma 5.5, $L(P)$ is a finite union of
languages of the form $N=\{u_0\}\{w_1\}^*\{u_1\}\ldots
\{w_k\}^*\{u_k\}$, with $u_j,w_j\in A^*$ and such that the rays
defined by $\eval w_1,\ldots,\eval w_k$ define an ordered simplex of
the triangulation $\cal T_{L_P}$ (recall that $L_P$ is the automatic
structure on $P$ determined up to equivalence by $L(P)$).  It follows
that the set $S$ of rays defined by the vertices of $C(P)$ is a subset
of $\partial L_P$. Since $\tau$ is conjugation invariant, replacing
the language $L$ by $g^{-1}Lg$ and restricting to $P$ gives the same
set $C(P)$ and hence the same set $S$.  However, $\partial L_P$ gets
replaced by $\partial L^g_P$.  Hence $S\subset\bigcap_{g\in G}\partial
L_P^g$.  Hence $\bigcap_{g\in G}\partial L_P^g$ is not contained in
any hemisphere.

We now show 3 implies 1. So suppose that for each maximal parabolic
subgroup $P$, $\bigcap_{g\in G}\partial L_P^g$ is not contained in any
hemisphere.  Let $P_1,\ldots,P_k$ be a set of representatives of the
conjugacy classes of maximal parabolic subgroups of $G$, with structure
$1\to\bbZ^{m_i}\to P_i\to F_i\to 1$.  We refer to the notation of the
proof of Theorem 4.2.  In particular, we start with any generating set
$B$ for $G$, which we are going to enlarge to a generating set $A$
that is appropriate for our purposes.

Denote $S_i:=\bigcap_{g\in G}\partial L_{P_i}^g$.  As in the proof of
Proposition 5.2 we can find an $F_i$-invariant polyhedron $Q_i\subset
\bbR^{m_i}$ with $0$ in its interior and with all its vertices
rational and on rays corresponding to points of $S_i$.  Choose a large
constant $K$ as in the proof of Theorem 4.2 and so that the set $A'_i$
of elements of $P_i$ that move the basepoint of $X$ at most distance
$K$ is a generating set for $P_i$.  Choose an integer $N_i$
sufficiently large so that the vertices of $N_i.Q_i$ have integer
coordinates and $A'_i\cap\bbZ^{m_i}\subset N_i.Q_i$. As in the proof
of Lemma 5.4, $A_i=A'_i\cup(N_i.Q_i\cap\bbZ^{m_i})$ is a good
generating set for $P_i$ with $C(A_i)=N_i.Q_i$.  Hence, by Lemmas 5.5
and 5.1, each language $L_{P_i}^g$ has an equivalent geodesic
representative in the generating set $A_i$.

We take the generating set $A=B\cup\bigcup_iA_i$ for $G$ as in the
proof of Theorem 4.2.  We claim that there exists a bound such that
elements of $L$ asynchronously fellow travel $A$-geodesics of $G$ at
distance given by this bound.  Indeed, if $u\in L$ and $v$ is a
geodesic $A$-word with the same value then, as in the proof of Theorem
4.2, we may write $u$ and $v$ as $u=u_0p_0u_1p_1\ldots$ and
$v=v_0q_0v_1q_1\ldots$ such that corresponding subwords $p_j$ and
$q_j$ lie in $A_{i_j}^*$ and, as portions of the paths $u$ and $v$,
begin and end a bounded distance apart, while corresponding portions
$u_j$ and $v_j$ of the paths $u$ and $v$ asynchronously
$\delta$-fellow travel for some $\delta$ that is independent of $u$.
The word $u_j$ lies in an automatic structure $L_j$ on $P_{i_j}$ which
depends only on the state of a machine for $L$ reached by the word
$u_0p_0\ldots p_{j-1}$.  Let $g_j=\eval{v_0q_0\ldots
q_{j-1}}^{-1}\eval{u_0p_0\ldots p_{j-1}}$. By what was said above, we
can replace each piece $q_j$ of $v$ by a word which is still geodesic
but lies in a language $L'_j$ equivalent to $g_jL_j{g_j}^{-1}$ on
$P_{i_j}$.  Since there are finitely many languages $L_j$ and $g_j$ is
of bounded size, there are finitely many languages $L'_j$ that need be
considered.  Thus $u$ asynchronously fellow travels the new word $v$
at distance bounded by the maximum of $\delta$ and the
fellow traveller constants between the $g_jL_jg_j^{-1}$ and $L'_j$.

Build an asynchronous 2-tape automaton $\cal T$ so that the language
of $\cal T$ is the set of pairs $(u,v)$ such that $u\in L$, $v$ is a
Cayley graph geodesic, $\eval u = \eval v$, and $u$ asynchronously
$K$-fellow travels $v$.  Consider the language $L_1$ which is the
projection onto the second factor.  It is regular (see, for example,
\cite{S1}). It has the asynchronous fellow traveller property and is
equivalent to $L$ since its words asynchronously fellow travel those
of $L$.  But it is a geodesic language, and a geodesic language with
the asynchronous fellow traveller property has the synchronous fellow
traveller property, so we are done.\qed

\rk Remark 5.7

Using the above ideas we can now easily give a proof of the result of
\cite{ECHLPT} that, if a biautomatic structure $L_i$ is chosen at each
cusp, there is a biautomatic structure on $G$ which restricts to
these.  Indeed, we take the generating set $A=B\cup\bigcup A_i$ for
$G$ of the above proof.  It has the falsification by fellow traveller
property so the language $N$ of geodesics is regular.  Moreover, by
Lemma 5.5 and Lemma 5.1, we can assume that $L_i$ is a geodesic
sublanguage of $A^*_i$. The desired biautomatic structure is then
$$\{w\in N \mid \text{ if $u$ is a maximal subword of $w$ with $u\in
A_i^*$ then $u\in L_i$}\}.$$
By the usual properties of regular languages, this is regular. The
asynchronous fellow traveller property follows from Lemma 3.2 and the
biautomaticity of the $L_i$.  The synchronous fellow traveller
property then follows from the fact that the language is geodesic.
\endrk

It is now easy to complete the proof of Theorem 3.1 by showing that
the map
$$\Phi\colon\Autstruct(G)\to\prod_{P\in\Par}\Autstruct(P)$$
has dense image.

\pf Proof of Theorem 3.1 completed

We suppose that $P'_1, \ldots, P'_k$ is a finite list of maximal
parabolic subgroups, and that structures $[L'_1],\ldots,[L'_k]$ are
chosen for these.  We must exhibit an automatic structure for $G$
which induces these.  We start with a biautomatic structure $[L]$ and
assume we have a generating set $A$ as in Lemma 4.6 and that $L$ is
geodesic.  We will modify $L$ to give an automatic structure $L'$
which induces the required structures on our chosen parabolics.  For
each $j=1,\ldots,k$ let $P'_j=g_jP_{i(j)}g_j^{-1}$ represent $P'_j$ as
a conjugate of one of $P_1,\ldots,P_m$.  The set $H_j$ of geodesic
words $x$ with $\eval x\in g_jP_{i(j)}$ which do not end in any
$A_{i(j)}$-letters is finite by Lemma 4.6. Moreover, for such an $x$
there is an automatic structure $L^{x}_j\subset A_{i(j)}^*$ on $\eval{
A_{i(j)}^*}$ so that $xL^{x}_jx^{-1} \sim L'_j$. We choose $L^y_j$ for
one $y\in H_j$ and then for any other $x\in H_j$ we put
$L^x_j=uL^y_ju^{-1}$, where $u\in A^*_{i(j)}$ is a word representing
$\eval{x^{-1}y}$. In this way we assure that the languages
$xL^x_jx^{-1}\in[L'_j]$ synchronously fellow-travel each other.

We take $L'$ to be the following language.
$$\eqalign{L'= &\{w\in L\mid w \text{ has no initial segment in
$\cup_j H_j$}\}\cup\cr
&\bigcup_{j=1}^k \bigcup_{x\in H_j} \{xp'v\mid \exists xpv\in L \text{ with
$p$ maximal in $A_{i(j)}^*$, and $p'\in L^x_j$}\}.}$$ 
$L'$ is regular by the usual properties of regular languages, and
clearly surjects to $G$.  It also clearly restricts to the correct
structures on the chosen parabolics.  Thus we need only verify the
fellow-traveller property.

We consider words of the form $xp'v$ as above. Let $\pi'$ and $\nu$ be the
$X$-geodesics for $p'$ and $v$ respectively and let $\gamma$ be the
$X$-geodesic for $p'v$.  Combining Lemmas A2.2 and A6, we see that the
path $\pi'\nu$ fellow-travels $\gamma$.  It now follows that $xp'v$
fellow-travels its $X$-geodesic after modification on horospheres.
Clearly, the same is true of the words in $L\cap L'$.

It is now easy to see that $L'$ has the asynchronous fellow-traveller
property.  For each word of $L'$ fellow-travels its $X$-geodesic after
modification on horospheres and the words of $L'$ all choose
equivalent structures on the horospheres they visit (we are using the
fact that $L$ is biautomatic). To see that $L'$ has the synchronous
fellow-traveller property it suffices to note that two nearby paths
spend a similar amount of time near any of the horospheres where we
have modified the structure.\qed

\HH Appendix:  Geometry in $X$

In this appendix we describe the geometry of the space $X$ of section
3.  Before we do so we recall some basic facts about hyperbolic space
$\bbH^n$ that we will need later.

Let $C$ be a closed convex subset of $\bbH^n$.  There is a retraction
$\sigma_C$ of $\bbH^n$ to $C$ by mapping any point to its closest
point in $C$.

\th Lemma A1

{\bf1.}~~If $x$ is a point distant $d$ from the closed convex set $C$ 
then $\sigma_C$ shrinks the local metric at $x$ by at least $e^d/2$.

{\bf2.}~~If $S$ is a horosphere of $\bbH^n$ and $x,y\in S$ are
hyperbolic distance $d$ apart, then their euclidean distance within
$S$ is $2\sinh(d/2)$.
\endth

\pf Proof

1.~~In fact, we will show that the degree of shrink is at least
$\cosh(d)$.  This exceeds $e^d/2$.  We first consider the case that
$C$ is a geodesic.  The degree of shrink at $x$ varies from a minimum
in the direction ``parallel'' to $C$ to a maximum of $\infty$ in the
direction towards $C$.  Using equation (7.20.3) of \cite{Be} it is
easy to show that the minimum shrink is by a factor of exactly
$\cosh(d)$.  Now let $C$ be any closed convex set.  Let $x$ and $x'$
be two points a very small distance apart, compared to their distance
from $C$, and let $y=\sigma_C(x),y'=\sigma_C(x')$.  Since the line
segment from $y$ to $y'$ lies in $C$, the angle $xyy'$ cannot be
acute, since if it where, moving $y$ towards $y'$ on this segment would
decrease $d(x,y)$.  Similarly, the angle $x'y'y$ is not acute. It follows
that $y$ and $y'$ lie at or between the projections of $x$ and $x'$ to
the geodesic $\gamma$ through $y,y'$.  Hence the degree of shrink for
$C$ is at least that for $\gamma$ and the lemma follows.

{2.} This follows from Theorem 7.2.1 of \cite{Be}.\qed

\th Lemma A2

{\bf1.}~~If $B$ is a horoball of $\Hyp^n$ and $x\in\Hyp^n$ then
any two geodesics from $x$ to $B$ first meet $B$ in points euclidean
distance less than $2$ apart in $S=\partial B$.

{\bf2.}If $B$ and $B'$ are disjoint horoballs and $\gamma$ is a
geodesic segment from a point of $S=\partial B$ to $S'=\partial B'$,
then the set of points that are on geodesic segments from points of
$S$ to points of $S'$ lies in a radius $3$ neighborhood of $\gamma$.
\endth

\pf Proof

{1.} We use the upper half-space model and carry $S$ to the horizontal
plane at height $1$.  The bound in question clearly increases as $x$
moves farther from $B$ with a supremum of $2$ when $x$ reaches the
boundary of upper half space.

{2.} We again position $S$ as above in the upper half-space model.
Call the set of points in question $N$.  It is a rotationally
symmetric solid that reaches its maximal diameter where it meets $S$
and $S'$.  Its intersection with $S$ is a euclidean disk $D$ which
increases in size as $B'$ approaches $B$.  In the extreme case that
$B$ and $B'$ are tangent, $B'$ is then drawn as a ball of radius $1/2$
in the model.  A hyperbolic geodesic which is tangent to both $B$ and
$B'$ is drawn as a semicircle of radius $1$ which meets the boundary
$\partial\Hypbar^n$ of the model in ``the shadow of $B'$'', that is,
in a point distance $\epsilon<1/2$ from the point of contact of $B'$
with the boundary.  The disk $D$ thus has euclidean radius
$1+\epsilon$, which is slightly less than $3/2$, so it has euclidean
diameter less than $3$, which implies hyperbolic diameter less than
$3$ (actually less than $2.4$).  \qed

We now return to the discussion of the geometry of $X$.  We will use
some basic facts about $\CAT(0)$ geodesic metric spaces which can be
found in \cite{Ba} or \cite{AB}.

We first recall the situation.  We have a geometrically finite group
$G$ acting on $\bbH^n$.  $\CH(G)$ is the convex hull for $G$, that is
the smallest non-empty convex subset of $\Hyp^n$ on which $G$ acts.
Each maximal parabolic subgroup $P$ of $G$ fixes a point at infinity
of $\bbH^n$ and hence fixes any horoball in $\bbH^n$ centered at this
point.  We choose a $G$-equivariant disjoint system of such horoballs,
one for each maximal parabolic subgroup and denote the horoball
corresponding to $P$ by $B_P$.  Then $X$ is the space
$$
X = \CH(G) - \bigcup_{P\in\Par} \int(B_P)
$$
with the path metric, that is the metric given by lengths of paths,
computed using the standard hyperbolic riemannian metric.  $X$ is
complete and locally compact, so it is a geodesic metric space. $G$
acts on $X$ by isometries, with finite stabilizers, and with compact
quotient.

We call the boundary piece $S_P=X\cap\partial B_P$ that results
from removing $\int(B_P)$ a \Em{horosphere of} $X$. Two horospheres
$S_P$ and $S_{P'}$ are have the same image in $M$ if and only if $P$
and $P'$ are conjugate in $G$.

Since $\CH(G)$ is convex, the metric on $\CH(G)$ is the restriction of
the metric on $\bbH^n$.  Now each horosphere is isometric to a convex
subset of euclidean space of dimension $n-1$.  It thus follows that
each point of the interior of $X$ lies in a neighborhood of curvature
$-1$, and each point of a horosphere of $X$ lies in a neighborhood
where the curvature is bounded above by $0$.  This makes $X$ a locally
$\CAT(0)$ space.  Since it is simply connected and locally compact, it
is a global $\CAT(0)$ space.  In particular, it is a geodesic metric
space, geodesics are unique, and they vary continuously with choice of
endpoint.

For points $x,y\in X$ we will denote $d_X(x,y)$ their distance apart
in $X$ and $d_{\bbH^n}(x,y)$ their hyperbolic distance, that is
distance in $\bbH^n$.

\th Lemma A3

Any two points $x,y\in X$ satisfy $d_{\bbH^n}(x,y)>2\log(d_X(x,y))$.
\endth

\pf Proof

Let $\gamma$ be the hyperbolic geodesic $x$ to $y$. We can find an
$X$-path from $x$ to $y$ by replacing each piece of $\gamma$ of the
form $\mu=\gamma\cap B_P$ by a geodesic $\mu'$ on the corresponding
horosphere. By Lemma A1.2, $\len(\mu')=2\sinh(\len(\mu)/2)$.  Since
$\sinh$ is a convex function and $\sinh(t)\ge t$ for all $t$, it
follows that the $X$-path we have created has length at most
$2\sinh(\len(\gamma)/2)$.  This is less than $\exp(\len(\gamma)/2)$,
so $d_X(x,y)<\exp(d_{\bbH^n}(x,y)/2)$.  The Lemma follows.\qed

Tatsuoka \cite{T} calls a path in $X$ a ``glancing geodesic'' if it is
$C^1$-smooth as a curve in $\bbH^n$ and decomposes piecewise into
hyperbolic geodesics in the interior of $X$ and euclidean geodesics on
the horospheres of $X$.  It is easy to see that any path which is not
a glancing geodesic can be shortened.  Since $X$ is a geodesic metric
space, the geodesics of $X$ are exactly the glancing geodesics.

Given a hyperbolic geodesic with endpoints in $X$, we can replace any
portion which leaves $X$ by the euclidean geodesic on the horosphere
whose horoball it enters.  We call such a path a \Em{rough geodesic}.
Such a path is not an $X$-geodesic, but it is roughly an $X$-geodesic
as the following Proposition shows.  It is actually a bit easier to
work with rough geodesics than with $X$-geodesics.

\th Lemma A4

There is a constant $\lambda$ so that every rough geodesic is a
$(\lambda,0)$-quasigeodesic in $X$.
\endth

\pf Proof

Suppose that $\gamma$ is an $\Hyp^n$ geodesic with endpoints in $X$.
We form $X'$ by deleting from $\CH(G)$ those horoballs that $\gamma$
enters.  We let $\rho$, $\sigma$, and $\tau$ be the $X'$, $X$, and
rough geodesics for $\gamma$ respectively.  We then have that
$\len(\rho)\le\len(\sigma)\le\len(\tau)$.  But the horospheres of $X$
(and hence those of $X'$) are bounded away from each other.  Using
Lemma A2 one easily checks that there is a global constant $\lambda$
so that $\len(\tau)\le \lambda\len(\rho)$ (this $\lambda$ is described
more closely in the next proof).  In particular $\tau$ is a
$\lambda$-quasigeodesic. \qed

\th Lemma A5

For any quasi-isometry constants $(\lambda,\epsilon)$ there exists a
constant $l$ such that if $w$ is a $(\lambda,\epsilon)$-quasigeodesic
in $X$ from $x$ to $y$ and $\gamma$ is the rough geodesic of $X$ from
$x$ to $y$ then $w$ asynchronously $l$-fellow travels a path obtained
from $\gamma$ by possibly modifying $\gamma$ on the horospheres it
visits. 
\endth

\pf Proof

Let $B_\gamma$ be the union of the horoballs that $\gamma$ visits and
$C_\gamma$ be the convex hull in $\Hyp^n$ of the union $\gamma\cup
B_\gamma$.  By Lemma A2, $C_\gamma$ is contained in the union of
$B_\gamma$ and a $3$-neighborhood of $\gamma$. Let $r_\gamma$ be the
restriction $r_\gamma=\sigma_{C_\gamma}|_{X}$ of the retraction
$\sigma_{C_\gamma}$ of Lemma A1.

By Lemma A1, $r_\gamma$ shrinks the local metric at a point $x$ by at
least $e^d/2$, where $d$ is the hyperbolic distance of $x$ from
$C_\gamma$.  Thus, by Lemma A3, $r_\gamma$ shrinks the local metric at
$x$ by at least $d_X^2/2$, where $d_X$ is the $X$-distance of $x$ from
$C_\gamma$.

Now let $u$ be any path in $X$ which lies outside a $K$-neighborhood
of $C_\gamma$, and whose endpoints $a$ and $b$ are exactly distance
$K$ from $C_\gamma$. Denote $u_1$=$r_\gamma\circ u$.  Then we have
just shown that $\len(u_1)\le 2\len(u)/K^2$.  Now $C_\gamma$, and
hence also $u_1$, may not lie completely in $X$.  We can correct this
as follows.  Note that $u_1$ runs on the boundary of $C_\gamma$. It
may stray from $X$ while running between horoballs $B$ and $B'$ that
$\gamma$ visits.  Let $u_1'$ be the portion of $u_1$ that does this.
We may replace $u_1'$ by a curve that runs from the start of $u_1'$
along the horosphere $\partial B$ to $\gamma$, along $\gamma$ to
$\partial B'$ and then along $\partial B'$ to the endpoint of $u_1'$.
By Lemma A2, this increases the length of $u_1'$ by a factor $k$ that
only depends on the minimal distance $\beta$ between distinct
horoballs ($k=1+6/\beta$ suffices). Also, $u_1$ may also stray from $X$
when running from one of its end points to a horoball.  Again, it is
easy to see that the relevant portion can be modified to stay in $X$
while increasing its length by at most the above factor $k$.  Finally,
$u_1$ might depart a horoball $B$, stray from $X$, and then return to
$B$.  In this case the relevant portion can be replaced by a shorter
geodesic on $\partial B$.  Summarizing, $u_1$ can be replaced by a
path $u_2$ in $X$ with the same endpoints of length $\len(u_2)\le
k\len(u_1)\le 2k\len(u)/K^2$.  By connecting the endpoints $a,b$ of
$u$ with the endpoints of $u_2$, we replace $u$ by a path $u'$ in $X$
of length at most $2(K+k\len(u)/K^2)$.  A simple calculation now shows
that if $K$ exceeds both $\sqrt{2\lambda}$ and $\epsilon/2$ and
$\len(u)> K':=\lambda K^2(2K-\epsilon)/(K^2-2\lambda)$ then
$\len(u')<\len(u)/\lambda-\epsilon$, so $u$ is not
$(\lambda,\epsilon)$-quasigeodesic.  Choose such $K$ and $K'$.  Then,
with $l_0=K+K'/2$, if the path $w$ of the lemma strayed outside the
$l_0$-neighborhood of $C_\gamma$ it would stray outside the
$K$-neighborhood for length at least $K'$.  It would thus fail to be
$(\lambda,\epsilon)$-quasigeodesic.

Thus $w$ stays in an $l_0$-neighborhood of $C_\gamma$, so it stays in
an $l_0$-neighborhood of the path $r_\gamma\circ w$, and hence, by
lemma 2A, it stays in an $(l_0+3)$-neighborhood of some path
$\gamma_1$ obtained by modifying $\gamma$ to follow $r_\gamma\circ w$
on horospheres. 

Now it is a standard fact that if two quasigeodesics run in a bounded
neighborhood of each other then they asynchronously fellow travel at
a distance depending only on the size of the neighborhood and
quasigeodesity constant.  Thus the lemma is proved.\qed

Applying Lemma A5 to an $X$-geodesic we see

\th Lemma A6

There exists a constant $\delta$ such that every rough geodesic
asynchronously $\delta$-fellow-travels its $X$-geodesic.\qed
\endth

In fact, it seems pretty clear that $\delta=1.5$ suffices, though our
proof does not give this.
\medskip

It is worth describing triangles in $X$, although we do not use this
in this paper.  Our description is rather sharper than the one in
\cite{T}. We first recall a characterization of triangles in a
hyperbolic metric space. Suppose that $\Delta=\alpha\beta\gamma$ is
such a geodesic triangle.  Then there is $\delta$, depending only on
the space, such that the sides of $\Delta$ decompose as
$\alpha=\alpha_0h_\alpha\alpha_1^{-1}$,
$\beta=\beta_0h_\beta\beta_1^{-1}$,
$\gamma=\gamma_0h_\gamma\gamma_1^{-1}$ with the following properties:

\Item{\bull} $\len(\alpha_1)=\len(\beta_0)$,
$\len(\beta_1)=\len(\gamma_0)$, and $\len(\gamma_1)=\len(\alpha_0)$.

\Item{\bull} Each of the pairs $\alpha_1$, $\beta_0$; $\beta_1$,
$\gamma_0$; $\gamma_1$, $\alpha_0$ $\delta$-fellow travel.

\Item{\bull} The portion of $\Delta$ labelled $h_\alpha$, $h_\beta$,
$h_\gamma$ has diameter less than $\delta$.  

Thus, we can think of $\Delta$ as consisting of a bounded \Em{hub}
with three thin \Em{spokes}.  One says $\Delta$ is \Em{$\delta$-thin}.
Geodesic triangles in $X$ behave
similarly, except that here, the hub may approximate a large euclidean
triangle.  We leave the proof of the following to the reader.

\th Proposition

There is $\delta$ with the following properties.  Suppose that
$\Delta=\alpha\beta\gamma$ is a geodesic triangle in $X$. Then either
$\Delta$ is $\delta$-thin or the sides of $\Delta$ decompose as
$\alpha=\alpha_0h_{\alpha}\alpha_1^{-1}$,
$\beta=\beta_0h_{\beta}\beta_1^{-1}$,
$\gamma=\gamma_0h_{\gamma}\gamma_1^{-1}$ with the following
properties:

\Item{\bull} $\len(\alpha_1)=\len(\beta_0)$,
$\len(\beta_1)=\len(\gamma_0)$, and $\len(\gamma_1)=\len(\alpha_0)$.

\Item{\bull} Each of the pairs $\alpha_1$, $\beta_0$; $\beta_1$,
$\gamma_0$; $\gamma_1$, $\alpha_0$ $\delta$-fellow travel.

\Item{\bull} The portion of $\Delta$ labelled $h_{\alpha}$, $h_{\beta}$,
$h_{\gamma}$ lies on a horosphere of $X$ and may be extended to a
geodesic euclidean hexagon $h_{\alpha}c_1
h_{\beta}c_2h_{\gamma}c_3$ where each of the $c_i$ has length at
most $\delta$. \qed\endth

In fact, if $\Delta'$ is the triangle in $\Hyp^n$ with the same
vertices as $\Delta$ then the second case of the Proposition occurs
when the hub of $\Delta'$ lies entirely in a horoball of $X$.

\Bib{References} 
 
\MaxReferenceTag{ECHLPT}

\rf{AB} S.B.~Alexander and R.L.~Bishop, The Hadamard-Cartan theorem in locally
convex spaces, l'Enseign.\  Math.~{\bf 36} (1990), 309--320.

\rf{ABC+}  J.M.~Alonso, T.~Brady, D.~Cooper,  V.~Ferlini, M.~ Lustig,
M.~Mihalik, M.~Shapiro, and H.~Short, H.~Short, ed., {\it Notes on
word hyperbolic groups}, in ``Group Theory From a Geometric
Viewpoint,'' E.~Ghys, A.~Haefliger and A.~Verjovsky, eds., World
Scientific (1990).

\rf{BGSS} G.~Baumslag, S.M.~Gersten, M.~Shapiro and H.~Short, 
Automatic groups and amalgams,  Journal of Pure and
Applied Algebra {\bf 76} (1991), 229--316. 

\rf{Ba} W.~Ballman, Chapter 10 of ``Sur les Groupes Hyperboliques
d'apr\`es Mickhael Gromov,'' E.~Ghys and P.~de la Harpe,eds., Progress
in Mathematics {\bf 83} Birkhauser Verlag (1990).

\rf{Be} A.F.~Beardon, The Geometry of Discrete Groups, Graduate Texts
in Mathematics 91, Springer-Verlag (1983).

\rf{Bo} B.~Bowditch, Geometric finiteness for hyperbolic groups,
Warwick PhD thesis (1988).

\rf{C} J.~Cannon, The combinatorial structure of cocompact discrete
hyperbolic groups, Geom.\  Dedicata {\bf16} (1984), 123--148.

\rf{ECHLPT}   D.B.A.~Epstein, J.W.~Cannon,   D.F.~Holt, S.V.F.~Levy
 M.S.~Paterson, and W.P.~Thurston, ``Word Processing in Groups,''
Jones and Bartlett Publishers, Boston (1992).

\rf{GS} S.M.~ Gersten and H.~Short, Rational sub-groups of 
bi-automatic groups, Annals of Math.\ {\bf 134} (1991), 125--158. 

\rf{M} J.~Milnor, A note on curvature and fundamental group, J.~Diff.\ 
Geom.\ {\bf 2} (1968), 1--7.

\rf{NS1} W.D.~Neumann and M.~Shapiro, Equivalent automatic structures
and their boundaries, Internat.\ J.~Alg.~Comp.\ {\bf 2} (1992), 443--469

\rf{NS2} W.D.~Neumann and M.~Shapiro, Automatic structures and
boundaries for graphs of groups, Internat.\ J.~Alg.~Comp.\ (to appear).

\rf{R} J.~Ratcliffe, The Foundations of Hyperbolic Manifolds, in
preparation.

\rf{Re} L. Reeves, Rational subgroups of geometrically finite
hyperbolic groups, preprint, University of Melbourne, 1993.

\rf{S1} M.~Shapiro, Non-deterministic and deterministic asynchronous
automatic structures, Internat.\ J.~Alg.\ and Comp.\ {\bf 3} (1992),
297--305.

\rf{S2} M.~Shapiro, Automatic structures and graphs of groups, in:
``Topology `90, Proceedings of the Research Semester in Low
Dimensional Topology at Ohio State,'' Walter de Gruyter Verlag,
Berlin - New York (1992), 355--380.

\rf{T} K.~Tatsuoka,  Finite volume hyperbolic groups are 
automatic, preprint (1990).

\endBib

\Coordinates 
The Ohio State University\\ 
Department of Mathematics\\ 
Columbus, OH 43210\\
USA
\par\noindent
The University of Melbourne\\
Department of Mathematics\\
Parkville, Victoria 3052\\
Australia
\endCoordinates 

\Coordinates 
City College\\ 
Department of Mathematics\\ 
New York, NY 10031\\
USA
\par\noindent
The University of Melbourne\\
Department of Mathematics\\
Parkville, Victoria 3052\\
Australia
\endCoordinates

\bye